\documentclass{amsart}
\usepackage{amsthm,amsmath,amsfonts,amssymb,amscd,mathrsfs,graphics}
\usepackage{latexsym}
\usepackage{graphicx}
\usepackage{rotating}

\usepackage{enumerate}

\newtheorem{thm}{Theorem}[section]

\newtheorem{lem}[thm]{Lemma}
\newtheorem{prop}[thm]{Proposition}
\newtheorem{cor}[thm]{Corollary}

\makeatletter
\renewcommand{\@seccntformat}[1]{\S{\csname
the#1\endcsname}\hspace{0.5em}}
\makeatother

\begin{document}

\title [Commutative Schur rings over $S_6$  ] {Commutative Schur rings over Symmetric groups II: the case $n=6$}

\author{ Amanda E. Francis, Stephen  P.  Humphries}
  \address{Department of Mathematics,  Brigham Young University, Provo, 
UT 84602, U.S.A.
E-mail: amandaefrancis@gmail.com, steve@mathematics.byu.edu} 
\date{}
\maketitle

\begin{abstract}  

We determine    the commutative Schur rings over $S_6$ that contain the sum of all the transpositions in $S_6$. There are eight such types (up to conjugacy), of which four have the set of all the transpositions as a principal set of the Schur ring.

\medskip

\noindent {\bf Keywords}: Schur ring, conjugacy class, finite group, group algebra.  \newline 
\medskip
\noindent {\bf AMS Classification}: Primary: 20C05; Secondary: 20F55
\end{abstract}

\section{Introduction}

For  a finite group $G$ and $X=\{x_1,\dots,x_k\} \subseteq G, |X|=k,$ we let $\overline X=x_1+\cdots +x_k\in \mathbb CG$. We also  let $X^{-1}=\{x^{-1}:x \in X\}$.
 
Let $\mathcal C_2$ denote the class of transpositions in the symmetric group $S_n$.
% We will say that  $X \subseteq \mathcal C_2$ is a {\it complete subset} if  there is  $Y \subseteq \{1,\dots,n\}$ such that $X=\{(i,j):i,j \in Y, i \ne j\}$. 
% The {\it size} of the complete subset is $|Y|$. 
 As a consequence of the main result of  \cite {hum}  we have:
 
%               \begin{thm} \label{thm992} Let $n \ge 6$ and   let 
% $\mathfrak S$ be a commutative  Schur-ring over $S_n$ that contains $\overline {\mathcal C_2}$. 
% Suppose that $C_1,\dots,C_h$ are the  principal sets of $\mathfrak S$ that are contained in $\mathcal C_2$. Then $h\le 3$, and  we may order the $C_i$ so that
%  $C_1,\dots,C_{h-1}$ are unions of disjoint complete subsets, and    $C_h=\mathcal C_2 \setminus \left ( \bigcup_{i=1}^{h-1} C_i\right)$. 
 
 % More precisely, for $h=1,2,3$ and $n \ge 6$, we have the following possibilities:

  %   \noindent(A) If $h=3$, then (up to a permutation of $\{1,\dots,n\}$) we have 
  %   $$C_1=\{(1,2)\}, \quad C_2=\{(i,j):3\le i<j\le n\}, \quad C_3=\mathcal C_2 \setminus (C_1 \cup C_2).$$

  %   \noindent (B) If $h=2$, then we may write $\mathcal C_2=C_1 \cup C_2$ where we have  one of the following:
 
%  (i) $n=6$ and $C_1$ is a union of two complete subsets of size $3$;

%  (ii) $n=6$ and $C_1$ is a union of three complete subsets of size $2$;   

 %  (iii)  $C_1$ is a   complete subset of size $n-1$.
  
  %   \noindent (C) If $h=1$, then  $C_1=\mathcal C_2$.
    
   %  Further, each of the above situations can be realized by a Schur-ring over $S_n$. 

 %  \end{thm}

% This gives:

\begin {cor} \label {corthm1} If
 $\mathfrak S$ is a commutative Schur ring over $S_6$ containing $\overline {\mathcal C_2}$, then $\mathfrak S$ determines  (up to conjugacy) one of the following partitions of $\mathcal C_2$:
 
 (i) $\mathcal C_2=\mathcal C_2$;
 
 (ii)  $\mathcal C_2=C_1 \cup C_2\cup C_3$ where 
 $$C_1=\{(1,2)\},\,\,\, C_2=\{(i,j):3\le i<j\le 6\} \text {  and }
C_3=\mathcal C_2 \setminus (C_1\cup C_2);$$   
 
 (iii) $\mathcal C_2=C_1 \cup C_2$ where 
 $$C_1=\{(1,2),(3,4),(5,6)\} \text { and }
C_2=\mathcal C_2 \setminus C_1;$$

 (iv) $\mathcal C_2=C_1 \cup C_2$ where 
 $$C_1=\{(1,2),(1,3),(2,3),(4,5),(4,6),(5,6)\}\text { and }
C_2=\mathcal C_2 \setminus C_1;$$

 (v) $\mathcal C_2=C_1 \cup C_2$ where 
 $$C_1=\{(i,j):1\le i<j\le 5\}\text { and } 
C_2=\mathcal C_2 \setminus C_1. $$\medskip
\end{cor} 
 
%  Let $\mathcal C_{2,2,2}\subset S_6$ denote the class of $(1,2)(3,4)(5,6)$. Then it is well-known (see e.g. \cite {km} \S 5.3) that there is an outer automorphism $\alpha$ of $S_6$ such that 
% $\alpha(\mathcal C_2) =\mathcal C_{2,2,2}$. 
% Thus the above options (i)-(v) of Corollary \ref {corthm1} also determines all the possible partitions for the class $\mathcal C_{2,2,2}$ in such a Schur-ring. 
% \medskip

If $H \le G$, then the set of orbits of elements of $G$ under the action of conjugation by elements of $H$ gives a Schur-ring that we denote $\mathfrak S(G,H)$.

In \cite {hum} we found all commutative Schur rings over $S_n, n\le 5$, that contain $\overline {\mathcal C_2}$. 
In this paper we  do the same for $S_6,$ this being our main result:

\begin{thm}\label{lem2222} The only commutative Schur rings over $S_6$ containing $\overline{\mathcal C_2}$ are (up to conjugacy):

(1) $Z(\mathbb CS_6)$;
%11

(2) $\mathfrak S(S_6,H_{120})$, where $H_{120}=\langle  (1, 4)(3, 5),
(1, 4, 6, 2, 5, 3)\rangle  \cong S_5$ is a $3$-transitive subgroup;
%19

(3) $\mathfrak S(S_6,H)$, where $H=A_6$;
%12

(4) $\mathfrak S_{36}$ (to be constructed in $\S 4$).
%12

(5)  $\mathfrak S(S_6,H)$, where $H=S_2 \times S_4\le S_6$;
%34

(6) $\mathfrak S(S_6,H)$, where $H=S_3 \wr  S_2\le S_6$;
%26

 (7)  $\mathfrak S(S_6,H)$, where $H=S_2 \wr  S_3\le S_6$;
 
(8)   $\mathfrak S(S_6,H)$, where $H=S_5$.

Of these, (1)-(4) are those where ${\mathcal C_2}$ is a principal set of the Schur ring.

 \end{thm}
 
 These Schur rings have dimensions $11,
 19,
 12,
 12,
    34,26,
  34,
 19$, respectively

We have shown in \cite  {hum2} that the dimension of a commutative  Schur-ring over $G$ is bounded by $s_G:=\sum_{i=1}^r d_i$, where the $d_i, i\le r,$ are the irreducible character degrees of $G$. In \cite {hum2} we showed that some (specific) groups realize this bound, while others do not. For example, we showed that $S_3,S_4,S_5$ each have a commutative Schur-ring  of this maximal dimension. We note that  $s_{S_6}=76$.
Then  Theorem \ref {lem2222} allows us to prove

\begin{cor}\label{corS66}
There is no commutative Schur-ring over $S_6$ of  dimension $76$.
 \end{cor}

 There is of necessity a certain amount of computation involved in the proof of Theorem \ref {lem2222}, however we have tried to restrict such computations to: 
 the Schur-ring algorithm that is described in $\S 2$; computations of Gr\"obner bases for certain ideals of polynomial algebras; a small number of small enumerations of possibilities.
 
 All computer computations involved in the preparation of this paper were accomplished using Magma \cite {ma}. 
 
\section {Schur-rings}

We now define the concept of a Schur ring  \cite {sch,wie,cur,mp}:

 A {\it Schur-ring} (or {\it S-ring}) over a finite group $G$ is a
sub-ring $\mathfrak S$ of $\mathbb{C}G$ that is constructed from a partition
$\{\Gamma _{1},\Gamma_2,\dots,\Gamma_m\}$ of the elements of $G$:
$
G=\Gamma_{1}\cup\Gamma_{2}\cup\dots\cup\Gamma_{m},$ with
$\Gamma_{1}=\{id\}$,  satisfying:

\noindent (1) if $1\le i\le m$,   then there is some $j \ge 1$ such that
$\Gamma_{i}^{-1} =\Gamma_{j}$;

\noindent (2)
 if  $1\le i,j\le m$, then
$\overline{\Gamma}_{i}\overline{\Gamma}_{j}=\sum_{k=1}^m
\lambda_{ijk}\overline{\Gamma }_{k},$ where $\lambda_{ijk} \in \mathbb Z^{\ge 0}$  for all $i,j,k$.

The $\Gamma_i$ are called the {\it principal sets} of the
$\text{S}$-ring.

\medskip

An S-ring naturally  gives rise to a subalgebra of $\mathbb CG$ by extending coefficients. We will usually think of S-rings as $\mathbb C$-algebras in this way.

We refer to the survey \cite {mp} for an account of recent developments and applications of the theory of Schur rings.

We recall the following fact (called the {\it Schur-Wielandt principle}, see Proposition 22.3 of \cite {wie}):

\begin{lem} \label{lemswie} 
Let  $\mathfrak S$ be an S-ring  over a group $G$. If $C \subseteq G$ satisfies $\overline{C} \in \mathfrak S$ and  $\sum_{g \in G}\lambda_g g \in \mathfrak S$, then 
   for all $\lambda \in \mathbb R$ the element $\sum _{g \in C} \delta_{\lambda_g,\lambda} g$ is in $\mathfrak  S$; here $\delta$ is the Kronecker delta function i.e. $\delta_{x,y}=0$ if $x \ne y$ and is $1$ otherwise.\qed 
   \end{lem}

          \noindent{\bf S-ring Algorithm}   Suppose that we have a subalgebra $H$ of $\mathbb C G$ ($|G|<\infty$) and we wish to find the smallest  S-ring that contains $H$.  Suppose that we start with a (ring) generating set 
            $c_1,\dots,c_r$ for the subalgebra  $H \subset \mathbb CG$.  For each $c_i$ partition the elements of $G$ according to their coefficients
    in $c_i$. For each such subset $C$ of this partition add $\overline{C}$  and $\overline{C^{-1}}$  to your   set of generators; do this for each $i$ and consider this  new set of  generators,  that we will denote by $d_1,\dots,d_t$. We simplify the set $d_1,\dots, d_t$ by eliminating any $\mathbb C$-linear dependences.  Now the $d_i$ determine and are determined by subsets of $G$, so that $t\le |G|$. Now consider the products $d_id_j, 1\le i,j \le t$. Again we partition the terms of $d_id_j$ according to their coefficients, and add in $\overline{C}, \overline{C^{-1}}$, for all sets $C$ in the partition. 
   Simplify this generating set using any linear dependences. This describes the basic step. If, after this basic step, one has a $\mathbb C$-basis for an S-ring, then we are done; otherwise we repeat the basic step. The process is guaranteed to terminate since $|G|<\infty$. It is easy to write a program to implement this algorithm (in, say, Magma \cite {ma}).

      \section {Covers of complete graphs}

 Let $K_n$ be the complete graph on $n$ vertices. Then for $\lambda \in \mathbb N$ and a graph $P$,  a {\it cover of $\lambda K_n$ by 
$P$'s} is a set $\mathfrak  T$ of subgraphs $P_1,\dots,P_m$ of $K_n$, each of which is isomorphic to $P$, and such that every edge of $K_n$ occurs $\lambda $ times in $P_1,\dots,P_m$. We will need the following result in $\S 4$. 

\begin{lem} \label{lemcovs63} 
Any cover $\mathfrak  T$ of $\lambda K_6$ by $r$ distinct triangles must have $5|r$. Further,
the cases $r=5,15$ do not happen, and the case $r=10$ is unique up to conjugation by an element of $S_6$:
\begin{align*}
%&\mathfrak C=
&\{     \{  1, 2, 3 \} ,
   \{  1, 2, 4 \} ,
     \{  2, 4, 5 \} ,
    \{  1, 3, 5 \} ,
     \{  1, 5, 6 \} ,    
 \{  3, 4, 5 \} , \{  3, 4, 6 \} ,
        \{  2, 5, 6 \} ,
             \\&\qquad \qquad  \{  2, 3, 6 \} ,   
          \{  3, 4, 6 \} ,  \{  1, 4, 6 \} 
      \}. \end{align*}
The situation $r=20$ is where we have all the triangles in $K_6$. 

\end{lem}\noindent{\it Proof} Let $\mathfrak T$ be a set of $r$ triangles giving the cover of $\lambda K_6$. Then counting edges we see that $3r=15\lambda$, so that $5|r$. 
Further, the number of triangles in $K_6$ is $20$.  If $r=15$, then taking the complement of $\mathfrak T$ in the set of all triangles gives the case $r=5$. Thus we consider the cases $r \in \{5,10\}$, since if $r=20$, then we have all the triangles 

\noindent {\bf Case A:}    $r=5$. Here $\lambda=1$. Then (by conjugating if necessary)  we can assume that $\{1,2,3\},\{1,4,5\} \in \mathfrak T$. Considering the edge $\{2,4\}$ one is forced (since  $\lambda=1$) to have $\{2,4,6\} \in \mathfrak T$. But now one checks that the edge $\{2,5\}$ cannot be in any triangle. Thus this case does not arise.\medskip

\noindent {\bf Case B:} $r=10$. Here $\lambda=2$. We consider two cases:

{\bf Case B1:} there is some $K_4\subset K_6$, all of whose triangles are in $\mathfrak T$. Here we can assume that the vertices of the $K_4$ are $1,2,3,4$. Then one sees that 
the only possibilities for extra triangles in $\mathfrak T$ are $\{i,5,6\}, i=1,2,3,4$ (otherwise we have edges with multiplicities greater than $2$).  But this makes it impossible to cover all edges
and have the multiplicity of the edge $\{5,6\}$ be $2$.\medskip

 Now assume that Case B1 does not happen. 

{\bf Case B2:} there is some $K_4\subset K_6$, all but one of whose triangles are in $\mathfrak T$. So assume that $\{1,2,3\},\{1,2,4\},\{1,3,4\} \in \mathfrak T$. 
Then, considering the edge $\{2,4\}$ we see that, up to permuting $5,6$, this forces $\{2,4,5\} \in \mathfrak T$.  
Considering the edge $\{2,5\}$  forces either (a) $\{2,5,6\} \in \mathfrak T$ or (b) $\{2,3,5\} \in \mathfrak T$. 

Assume that we have (a) $\{2,5,6\} \in \mathfrak T$.
Considering the edge  $\{1,6\}$   forces $\{2,3,6\} \in \mathfrak T$. 
Considering the edge  $\{2,6\}$   forces $\{1,5,6\} \in \mathfrak T$, and that this is the only triangle that can be in $\mathfrak T$ that contains $\{1,6\}$. 
Thus this case cannot happen.

If we have (b) $\{2,3,5\} \in \mathfrak T$, then considering the edge $\{2,6\}$ gives a contradiction.

 Now assume that Case B1 and Case B2 do not happen. Then we can assume that $\{1,2,3\},$ $ \{1,2,4\} \in \mathfrak T$, $\{1,3,4\},\{2,3,4\} \notin \mathfrak T$. 
Considering the edge $\{2,4\}$  forces $\{2,4,5\} \in \mathfrak T$ (up to permuting $5,6$).
Considering the edge $\{1,5\}$  (and using the fact that Case B2 does not occur) we must have $\{1,3,5\},\{1,5,6\} \in \mathfrak T$. 
Considering the edge $\{3,4\}$  (and using the fact that Case B2 does not occur) we must have $\{3,4,5\},\{3,4,6\} \in \mathfrak T$. 
Similarly, considering the edge $\{5,6\}$   we must have $\{2,5,6\} \in \mathfrak T$. 
Considering the edge $\{2,6\}$  we must have $\{2,3,6\} \in \mathfrak T$. This leaves $\{1,4,6\}$ as the remaining triangle.\qed 

\medskip

To each element  $(i,j,k) \in \mathcal C_3$ there is associated the triangle $\{i,j,k\}$  in $K_6.$ We let the following subset of $\mathcal C_3$  represent 
 the set of triangles in the above Lemma: 
\begin{align*}
&\mathfrak C_3=\{     (  1, 2, 3 ) ,
   (  1, 2, 4 ) ,
     (  2, 4, 5 ) ,
    (  1, 3, 5 ) ,
     (  1, 5, 6 ) ,    
 (  3, 4, 5 ) , (  3, 4, 6 ) ,
        (  2, 5, 6 ) ,
             \\&\qquad \qquad  (  2, 3, 6 ) ,   
          (  3, 4, 6 ) ,  (  1, 4, 6 ) 
      \}. \end{align*}

The set of triangles corresponding to the elements of $\mathfrak C_3$ can be  represented as the set of triangles of the hemi-icosahedron $\mathfrak H$ (a polyhedral decomposition of the projective plane). See Figure 1, where we have drawn the ten triangles ($2$-simplices)  of $\mathfrak H$, and the outside edges of these $2$-simplices are identified in pairs as usual.

   \begin{figure}[ht]
 \scalebox{.7}
{\includegraphics{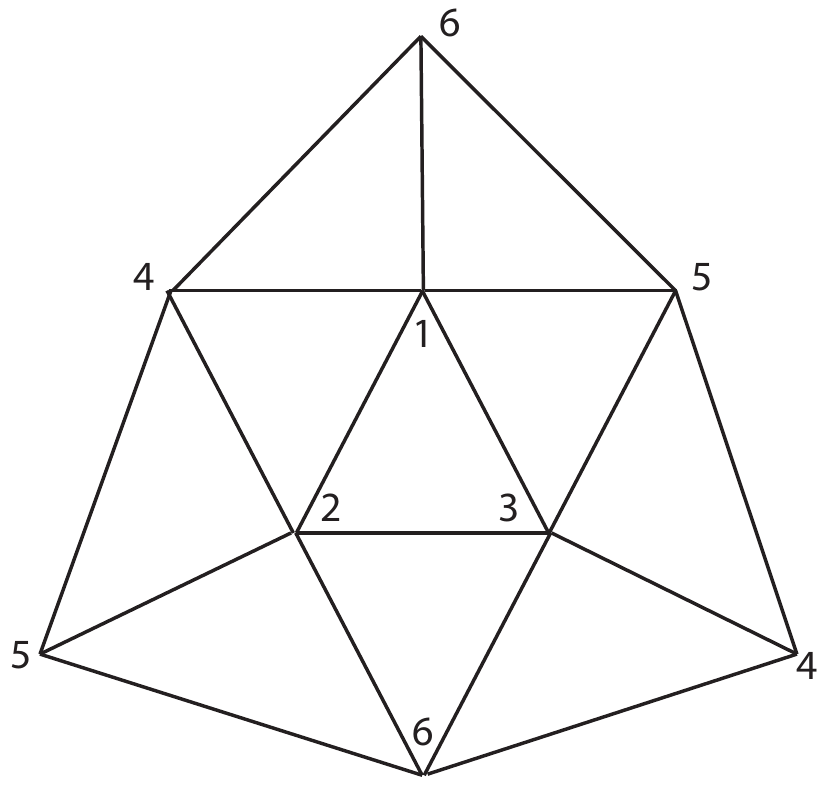}}
 \caption{ The   hemi-icosahedron corresponding to  $\mathfrak C_3$ } 
 \end{figure}
 
 The automorphism group of this $2$-complex is $A_6 \le S_6$, and it acts transitively on the ten triangles corresponding to the elements of  $\mathfrak C_3$.
\medskip

\section{ The case where $\mathcal C_2$ is a principal set.}
 
 In this section we assume that $\mathcal C_2$ is a principal set of a commutative S-ring $\mathfrak S$. This is the most difficult case.  Since $\overline {\mathcal C_2} \in \mathfrak S$
 it follows from Lemma 4.1 of  \cite {hum2} that $Z(\mathbb CG)$ is a subring of $\mathfrak S$.

 It follows from this and the Schur-Wielandt principle
 that if $C$ is a principal set of $\mathfrak S$, then $C$ is contained in some conjugacy class, namely the class of one of its elements.
 
 Let $\mathcal C_3$ denote the class of $3$-cycles in $S_n$.  In general for $\mu \vdash n$ we let $\mathcal C_\mu$ denote the class of elements of $S_n$ of cycle type $\mu$.
 
 \medskip 
 
 \noindent{\bf The class $\mathcal C_{3}$}.
 
 We first consider the principal sets  $C \subseteq  \mathcal C_3$ of $\mathfrak S$, with the goal of showing   $C=\mathcal C_3$. 
 
 For $\mu\vdash n$ and $\alpha=\sum_{g\in G} \lambda_g g$ we let 
$$\alpha_\mu=\sum_{g\in \mathcal C_\mu} \lambda_g g.$$

\begin{lem} \label{lemcov}  If $C \subseteq \mathcal C_3\subset S_n$ is a principal set  of $\mathfrak S$ and $\left (\overline{C}\cdot \overline{\mathcal C_2}\right )_{(2,1^{n-2})}
=\lambda \overline{\mathcal C_2}$, then $C$ determines a cover of $\lambda K_n$  by triangles.
Moreover, we have
$$3\cdot |C|=\lambda\cdot |\mathcal C_2|.$$
\end{lem}
\noindent {\it Proof} We have
$$(i,j,k) \cdot ((i,j)+(j,k)+(i,k))=(i,j)+(j,k)+(i,k).$$
Thus each $3$-cycle $(i,j,k)\in C$ contributes $(i,j)+(j,k)+(i,k)$ to the product    $\overline{C}\cdot \mathcal C_2$. 

Further,  for each $\alpha=(i,j,k)\in \mathcal C_3,$ there are only three $\beta \in \mathcal C_2$ (namely $\beta \in \{(i,j),(j,k),(k,i)\}$) with $\alpha\beta \in \mathcal C_2$.
Since each $(i,j)\in \mathcal C_2$ occurs $\lambda$ times in $\overline{C}\cdot \mathcal C_2$ we have a cover of $\lambda K_n$ by triangles, and   the first part of the result follows.
By counting edges we see that
$3\cdot |C|=\lambda\cdot |\mathcal C_2|.$
\qed\medskip

If $n=6$, then  Lemma \ref {lemcovs63} together with the fact that either $C=C^{-1}$ or $C \cap C^{-1}=\emptyset$, shows that (up to conjugacy) we have one of:

\noindent (I) $C=\mathcal C_3$;

\noindent (II) $|C|=20$ where $\overline {C}= \sum_{\alpha \in \mathfrak C_3} \alpha+ \alpha^{-1}$; 

\noindent (III) $|C|=10$ where $\overline {C}=\sum_{\alpha \in \mathfrak C_3} \alpha^{\varepsilon(\alpha)}$ with $\varepsilon(\alpha) \in \{1,-1\}$; 

\noindent (IV) $\overline {C}=\sum_{1\le i<j<k\le 6} (i,j,k)^{\varepsilon(i,j,k)}$, for some  $\varepsilon(i,j,k) \in \{1,-1\}$.
\medskip

We want to show that (I) is the only possibility. If we have (II), then applying the S-ring algorithm to the algebra generated by $Z(\mathbb CS_6)$ and $\overline {C}$ gives an S-ring that is not commutative. Thus (II) does not happen.
\medskip

Suppose that  we have (III), so that $\overline {C}=\sum_{\alpha \in \mathfrak C_3} \alpha^{\varepsilon(\alpha)}$ with $\varepsilon(\alpha) \in \{1,-1\}$. 

We will use Figure 1 to determine a {\it standard orientation} (counter clockwise) for each of the triangles of $\mathfrak C_3$, so that, for example the order $1,2,3$ determines a positive orientation for the triangle $\{1,2,3\}$ in $\mathfrak C_3$. 

Since we have (III), every $(i,j,k) \in C$ has the form $\alpha^{\varepsilon(\alpha)}$, for some $\alpha \in \mathfrak C_3$.
Each such 
 $\alpha^{\varepsilon(\alpha)}=(i,j,k)\in C$ determines a triangle $\mathfrak t(\alpha^{\varepsilon(\alpha)})=\{i,j,k\}$ of $\mathfrak H$. In fact $\alpha^{\varepsilon(\alpha)}=(i,j,k)\in C$ determines an {\it orientation} for  $\mathfrak t(\alpha^{\varepsilon(\alpha)})$, this being 
  $i\mapsto j  \mapsto k\mapsto i$. 
 % For $C \subset \mathcal C_3, |C|=10,$ as in (III), we see that each $\alpha \in C$ determines an orientation of a triangle of $\mathfrak H$. 
 We denote the set of oriented triangles of $\mathfrak H$ determined by $C$ as $\mathfrak H(C)$.  Further, each  oriented triangle determines an orientation on each of its edges.
  
   An edge of $\mathfrak H$ will be said to be {\it oppositely-oriented in $\mathfrak H(C)$} if that edge has different orientations in the two triangles of $\mathfrak H(C)$ in which  it is contained.
 Let $o(C)$ denote the number of oppositely-oriented edges in $\mathfrak H(C)$.  We note that $o(C)<15$, since  $\mathfrak H$ has $15$ edges.
 
 \begin{lem} \label {lem1}
 (a) If $(i,j,k), (r,s,t) \in \mathcal C_3$, then $(i,j,k)(r,s,t) \in \mathcal C_3$ if and only if either 
 
 \noindent (i) $(i,j,k)=(r,s,t)$; or 
 
 \noindent (ii) 
 $\{i,j,k\} \cap \{r,s,t\}=\{u,v\}$ has size $2$ and the edge $u,v$ in $(i,j,k)$ is oriented differently from that in 
 $(r,s,t).$
 
 (b)  Two (oriented) triangles of  $\mathfrak H(C)$ that share an edge will contribute two elements of $\mathcal C_3$ to $\overline{C}^2$ if and only if their common edge is oriented differently in the two triangles, if and only if the two triangles are oriented the same. In any other  case they contribute no elements of $\mathcal C_3$.
 
 (c)   $\left ( \overline {C}^2-\overline{C^{-1}}\right )_{(3,1^3)}$ is a sum of  $2o(C)$ elements of $\mathcal C_3$. 
 
 \end{lem} 
 \noindent{\it Proof} (a) We may assume that $i,j,k,r,s,t\le 6$, and one checks a small number of cases. 
 
(b)
 Now applying Lemma \ref {lem1} (a) to $\overline{C}^2|_{(3,1^3)}$ we see that there is the contribution that comes from  
 $\sum_{\alpha \in C} \alpha^2$. We note that none of these elements $\alpha^2$ of $\mathcal C_3$ are in $C$, since $C \cap C^{-1}=\emptyset$. 
 Any other product comes from a pair $\alpha_1,\alpha_2 \in C$ that satisfy Lemma \ref {lem1} (a) (ii). Thus $\alpha_1,\alpha_2$ determine an orientation of two of the triangles of $\mathfrak H$
 that have opposite orientations on the common edge of these two triangles.
 
 So considering $\mathfrak H(C)$ we see that two triangles of  $\mathfrak H(C)$ that share an edge will contribute two elements of $\mathcal C_3$ to $\overline{C}^2$ if and only if their common edge is oriented differently in the two triangles, if and only if the two triangles are oriented the same.

 (c) This follows from (b) and its proof.  \qed\medskip
 
 Let $\Omega(C)$ be the subgraph of the $1$-skeleton $\mathfrak H(C)^{(1)}$ of $\mathfrak H(C)$  consisting of edges that are oppositely-oriented.
 Let $f_+, f_-$ be the number of triangles  of $\mathfrak C_3$  that are oriented positively or negatively (respectively)  in $\mathfrak H(C)$. Then $f_++f_-=10$ and we may assume that $f_+\ge f_-$ (else replace $C$ by $C^{-1}$).
 
 \begin{lem} \label {lem2} (a) The graph $\Omega(C)$  contains all the vertices of $\mathfrak H(C)^{(1)}$. 
 
 (b) $o(C)\in \{5,10\}$.

 (c)  $f_+\equiv f_-\equiv o(C) \mod 2$. In particular, 
 
  if $f_-=0$, then $o(C)=12$; 
 
 if $f_-=1$, then $o(C)\in \{9,11\}$; 
 
 if $f_-=2$, then $o(C)\in \{6,8,10\}$;
 
  if $f_-=3$, then $o(C)\in \{5,7,9\}$;
  
    if $f_-=4$, then $o(C)\in \{6,8,10\}$;

     if $f_-=5$, then $o(C)\in \{3,7\}$;
     
   if $f_-=6$, then $o(C)=6$.

(d)  If $o(C)=5$, then  in Figure 2 we show the two possibilities for $C$ (up to a permutation of $1,\dots,6$) where we have shaded the   positively (say) oriented triangles, and where the oppositely oriented edges are drawn dashed.
 
    \begin{figure}[ht]
 \scalebox{.6}
{\includegraphics{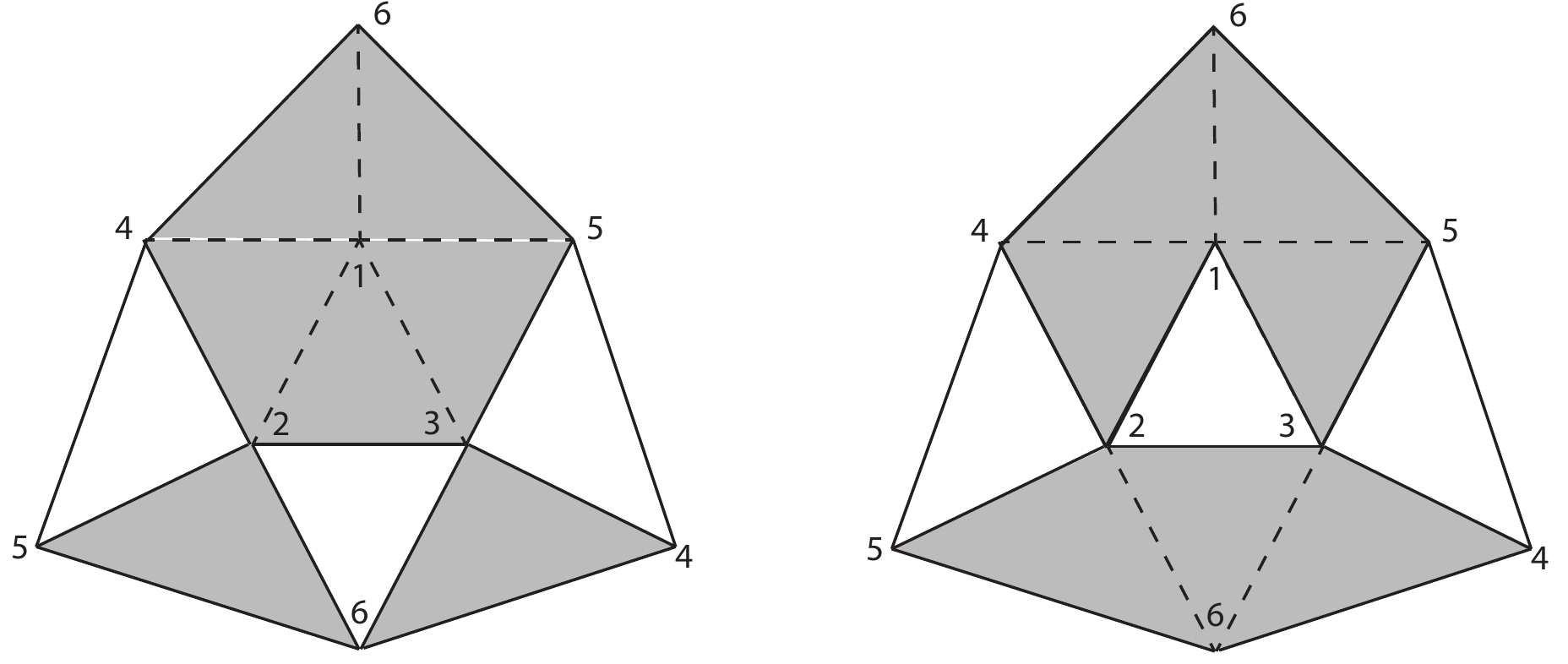}}
 \caption{ Orientations of $\mathfrak H(C)$ giving $o(C)=5$. } 
 \end{figure}

\end{lem} 
  \noindent{\it Proof} (a) This follows, since every edge of $\mathfrak H(C)^{(1)}$  has degree $5$, so that for every vertex $v$  of $\mathfrak H(C)$, at least two adjacent  triangles that share the vertex $v$ must have the same orientation in 
  $\mathfrak H(C)$. 
 
 (b) We know that every principal set  contained in
$\mathcal C_3$ has size a multiple of $5$.  Lemma \ref {lem1} (c) shows that $10|o(C)$. But $o(C)<15$
then gives the result. 
 
 (c) Changing the orientation of any triangle of $\mathfrak H$ changes the orientation of the three edges of that triangle. It follows that  $f_+\equiv f_-\equiv o(C) \mod 2$. The rest follows by looking at each particular case.

 (d) We use (c) to check this. \qed\medskip

 Thus Case (III) has been reduced to showing that the two S-rings determined by Lemma \ref {lem2} (d) do not give commutative S-rings. One checks this latter fact using the S-ring algorithm. This concludes consideration of (III).
 \medskip
 
 Now assume that we have (IV), so that any  subset of $\mathcal C_3$ that is a principal set has size $20$. Thus $C^{-1}\ne C$ is also a principal set. Thus we have    
 $$\overline {C}=\sum_{1\le i<j<k\le n} (i,j,k)^{\varepsilon(i,j,k)} \text { where } \varepsilon(i,j,k) \in \{1,-1\}.$$

\begin{lem} \label {lemuu}

( Let $a,b,c \in \mathcal C_3$. Then at most one of the eight products 
$a^{\pm 1}b^{\pm 1},  b^{\pm 1}a^{\pm 1}$ can be equal to $c$.

% (b) If $i<j<k\le n, r<s<t\le n, u<v<w\le n$ and $c=(i,j,k), a=(r,s,t), b=(u,v,w)$, then we have $ab=c$ exactly for the following cases:
% \begin{align*}
% &(i,j,k)=(j,k,m)(i,j,m), \text{ where } 1\le k<m\le n;\\
% &(i,j,k)=(i,j,m)(i,m,k), \text{ where } j <m<k;\\
 % &(i,j,k)=(i,m,k)(m,j,k), \text{ where } i<m<j;\\
 % &(i,j,k)=(m,j,k)(m,i,j), \text{ where } 1\le m<i.\\
 %  \end{align*}

% (c) If $i<j<k\le n, r<s<t\le n, u<v<w\le n$ and $c=(i,j,k), a=(r,s,t), b=(u,v,w)$, then we have $a^{-1}b=c$ exactly for the following cases:

% \begin{align*}
% &(i,j,k)=(i,k,m)^{-1}(j,k,m), \text{ where }  k<m\le n;\\
% &(i,j,k)=(j,m,k)^{-1}(i,j,m), \text{ where } j <m<k;\\
 % &(i,j,k)=(i,m,j)^{-1}(i,m,k), \text{ where } i<m<j;\\
 % &(i,j,k)=(m,i,k)^{-1}(m,j,k), \text{ where } 1\le m<i.\\
 %  \end{align*}

% (d) 
 % If $i<j<k\le n, r<s<t\le n, u<v<w\le n$ and $c=(i,j,k), a=(r,s,t), b=(u,v,w)$, then we have $ab^{-1}=c$ exactly for the following cases:

% \begin{align*}
% &(i,j,k)=(i,j,m)(i,k,m)^{-1}, \text{ where } 1\le k<m\le n;\\
% &(i,j,k)=(i,m,k)(j,m,k)^{-1}, \text{ where } j <m<k;\\
%  &(i,j,k)=(m,j,k)(i,m,j)^{-1}, \text{ where } i<m<j;\\
%  &(i,j,k)=(m,i,j)(m,i,k)^{-1}, \text{ where } 1\le m<i.\\
%   \end{align*}

% (e)  If $i<j<k\le n, r<s<t\le n, u<v<w\le n$ and $c=(i,j,k), a=(r,s,t), b=(u,v,w)$, then we have $a^{-1}b^{-1}=c$ exactly when $a=b=c$.

\end {lem}

\noindent {\it Proof } (a) We may assume that $c=(1,2,3)$, and then one only has to check the cases where $a,b \in \mathcal C_3 \cap S_4$.
\qed\medskip

For  $c=(i,j,k)\in C$ we  now wish to determine the coefficient of $c$ in $\overline {C}^2$ as a polynomial  function of the $\varepsilon(r,s,t)$. 

First assume that $\varepsilon(i,j,k)=1$. Suppose   $a=(r,s,t)^{\delta_1},$ $ b=(u,v,w)^{\delta_2}, $ are distinct where $\delta_1=\delta_{1,i,j,k,r,s,t,u,v,w},\delta_2=\delta_{2,i,j,k,r,s,t,u,v,w} \in \{1,-1\},  \, r<s<t, u<v<w,$ and $c$ is one of 
$ab,  ba$.  Then by Lemma \ref {lemuu} (a)  the coefficient of $c$ in $(a+b)^2$ is $1$. Thus if $a,b \in C$, then they will contribute $1$ to the coefficient of $c$ in 
   $\overline{C}^2$.
   So when  $a=(r,s,t)^{\varepsilon(r,s,t)}, b=(u,v,w)^{\varepsilon(u,v,w)} \in C$  are distinct, then  they will contribute $1$ to the coefficient of $c$ in 
   $\overline{C}^2$ if and only if $\varepsilon(r,s,t)=\delta_1, \varepsilon(u,v,w)=\delta_2$.

Thus we have (i) of

 \begin{lem} \label {lemrstuvw}  (i)  When   $c=(i,j,k)\in C$,   the coefficient of $c=(i,j,k)$ in    $\overline{C}^2$ coming from the pair $(r,s,t)^{\varepsilon(r,s,t)},(u,v,w)^{\varepsilon(u,v,w)}$ is
$$\frac {(1+\delta_1\varepsilon(r,s,t))}{2} \cdot \frac {(1+\delta_2\varepsilon(u,v,w))}{2}.$$

(ii) When   $c=(i,j,k)\notin C$,   the coefficient of $c$ in    $\overline{C}^2$ coming from the pair $(r,s,t)^{\varepsilon(r,s,t)},(u,v,w)^{\varepsilon(u,v,w}$ is
$$\frac {(1-\delta_1\varepsilon(r,s,t))}{2} \cdot \frac {(1-\delta_2\varepsilon(u,v,w))}{2}.$$

(iii)  For any   $c=(i,j,k)^{\varepsilon(i,j,k)}, 1\le i<j<k\le6$,   the coefficient of $c$ in    $\overline{C}^2$ coming from the pair $a=(r,s,t)^{\varepsilon(r,s,t)}, b=(u,v,w)^{\varepsilon(u,v,w}$ is
$$\frac {(1+\delta_1\varepsilon(r,s,t)\varepsilon(i,j,k))}{2} \cdot \frac {(1+\delta_2\varepsilon(u,v,w)\varepsilon(i,j,k))}{2}.$$
\end{lem}
\noindent {\it Proof} The proof of (ii) is similar to what we did above for (i), and (iii) is a restatement of  (i) and (ii).\qed\medskip 

Now let $P$ be a polynomial ring over $\mathbb Q$ with generators $X(i,j,k)$ where $1\le i<j<k\le 6$ and define, for $c=(i,j,k)^{\varepsilon(i,j,k)}, a=(r,s,t)^{\varepsilon(r,s,t)}, b=(u,v,w)^{\varepsilon(u,v,w)}, a \ne b,$ where $c$ is one of $ab, ba$.  
We define the polynomial 
$$S(a,b,c)=\frac {(1+\delta_1X(r,s,t)X(i,j,k))}{2} \cdot \frac {(1+\delta_2X(u,v,w)X(i,j,k))}{2}.$$
Here  $\delta_1=\delta_{1,i,j,k,r,s,t,u,v,w},\delta_2=\delta_{2,i,j,k,r,s,t,u,v,w} \in \{1,-1\}, \, r<s<t, u<v<w,$ as in the above. Thus the $X(i,j,k)$ correspond to the $\varepsilon(i,j,k)$. 
We let  $E_\varepsilon:P \to \mathbb Q$ be the evaluation map $X(i,j,k) \mapsto \varepsilon_{i,j,k}$.

For each $c=(i,j,k)^{\varepsilon(i,j,k)}, i<j<k,$ we let $S(c)$ be the sum of all $S(a,b,c)$ where  $a=(r,s,t)^{\varepsilon(r,s,t)}, b=(u,v,w)^{\varepsilon(u,v,w)} \in C$  are distinct and one of $ab, ba$ is $c$. Thus $S(c)$ is a sum of nine of the $S(a,b,c)$s. It is clear that if $c=(i,j,k)^{\varepsilon(i,j,k)}, i<j<k,$ and $c'=(i',j',k')^{\varepsilon(i',j',k')}, i'<j'<k',$ then $E_\varepsilon(S(c))=E_\varepsilon(S(c'))$, since $C$ is a principal set. 

We now define the ideal $I$ whose generators are all $S(c)-S(c')$, for $c,c' $ as above, and all $X(i,j,k)^2-1$ for 
$ 1\le i<j<k\le 6.$ Since we can always conjugate $C$ so that $\varepsilon(1,2,3)=1=\varepsilon(4,5,6)$ we also let $X(1,2,3)-1,X(4,5,6)-1$ be in $I$. 

Now if there is such a subset $C$, then the ideal $I$ will not be equal to $P$. However a 
Gr\"obner basis calculation (\cite {ma})  shows that $I=P$, and so (IV) is not possible. Thus we now have:

 \begin{prop} \label{prop23} Any commutative S-ring over $S_6$ that contains $\mathcal C_2$ as a principal set, also contains $\mathcal C_3$ as a principal set.
\qed \end{prop}

\noindent{\bf The class $\mathcal C_{(2^2,1^2)}$}.

We let $\mathcal C_{2,2}$ denote  $\mathcal C_{(2^2,1^2)}$.
 Here we  have:

\begin{lem} \label{lem22cycle}
If $C \subset \mathcal C_{2,2}$ satisfies $\overline{C} \in \mathfrak S$, where $\mathfrak S$ is a commutative S-ring over $S_n$ containing $\mathcal C_2$ as a principal set, then 
$C$ determines a cover of $\lambda K_n$ by pairs of disjoint edges.\end{lem}
\noindent {\it Proof}
This follows from
\begin{align*}
[(i,j)(k,m) \cdot \overline { \mathcal C_2}]_{(2,1^{n-2})}&=\left [(i,j)(k,m) \cdot [(i,j)+(k,m)+\cdots]\right ]_{(2,1^{n-2})} \\&= (i,j)+(k,m), 
\end{align*} together with the fact that $\mathcal C_2$ is a principal set, so that $[\overline{C}\cdot \overline{\mathcal C_2}]_{(2,1^{n-2})}=\lambda\overline{\mathcal C_2}.$\qed

\medskip 
 
 In the above situation we have: $2|C|=\lambda |\mathcal C_2|$. 
 
For each of $n=4, \lambda=1; n=5, \lambda=1,$ and $n=5,\lambda=2,$ one can check that  there is a unique such cover.
%$(1,2)(3,4)+ (1,3)(2,4)+ (1,4) (2,3).$
%For $n=5$ and $\lambda=1$ we have the unique (up to conjugacy)  cover
%\begin{align*}&(1,2)(3,4)+(1,3)(4,5)+(1,4)(2,5)+(1,5)(2,3)+(2,4)(3,5).\end{align*}
%For $n=5$ and $\lambda=2$ we have the unique (up to conjugacy) cover
%\begin{align*} &(1, 3)(2, 4)+(1, 5)(2, 4)+(1, 5)(2, 3)+(2, 3)(4, 5)+(1, 3)(4, 5)\\&+(2, 5)(3, 4)+(1, 2)(3, 4)+(1, 2)(3, 5)+(1, 4)(3, 5)+(1, 4)(2, 5).\end{align*}
For $n=6$ we have $2|C|=15\lambda$, so that  $\lambda$ must be even, and $15$ divides $|C|$.
 We note that $|\mathcal C_{(2^2,1^2)}|=45$, so that we can assume  $|C|<45$. 
 If $|C|=30$, then we replace $C$ by $\mathcal C_{(2^2,1^2)} \setminus C$, so that we can assume $|C|=15$.
We now show 
\begin{lem}\label{lem333} For $C \subset \mathcal C_{(2^2,1^2)}, |C|=15,$ as above, the only possibilities for $\overline {C}$ are those  conjugate to the following element (or its complement in $\mathcal C_{(2^2,1^2)}$)
\begin{align*}
& \,\,\,\,\,\,\,
(1, 2)(3, 4)+(1, 2)(5, 6)+(3, 4)(5, 6)\\&+
 (1, 3)(2, 5)+(1, 3)(4, 6)+(2, 5)(4, 6)\\&+
 (1, 5)(2, 4)+(1, 5)(3, 6)+(2, 4)(3, 6)\\&+ 
 (1, 6)(2, 3)+ (1, 6)(4, 5)+ (2, 3)(4, 5)\\&+
 (1, 4)(2, 6)+(1, 4)(3, 5)+(2, 6)(3, 5).\end{align*}

The stabilizer of this element (under  conjugation by elements of $S_6$) is $$H_{120}=\langle (1, 4)(3, 5),
(1, 4, 6, 2, 5, 3)\rangle,$$ so that $[S_6:H_{120}]=6$, and there are $6$ conjugates of this element.

Further, 
the S-ring generated by $Z(\mathbb CG)$ and the above element is $ \mathfrak S(S_6,H_{120})$.
 
\end{lem}
\noindent{\it Proof} The last two paragraphs are simple computations once we have the first part.

For the first part, let $K_n^e$ denote the graph with vertices the edges $\{i,j\}$ of $K_n$, and where we have  an edge of $K_n^e$ whenever the corresponding vertices (edges of $K_n$)  are disjoint.
Thus each edge $\{i,j\} - \{k,m\}$ of $K_n^e$ corresponds to an element $(i,j)(k,m) \in \mathcal C_{(2^2,1^{n-4})}$, so that we can refer to an edge of $K_n^e$ by that element.

Now let $C \in \mathcal C_{(2^2,1^{2})} $ determine a  cover of $2K_6$ by subgraphs isomorphic to $K_2 \cup K_2$. Then every vertex of $K_6^e$ is in exactly two edges. Thus
$C$ determines a set of disjoint cycles $\Gamma_1,\dots,\Gamma_r$ in $K_6^e$. Each cycle $\Gamma_i, i \le r,$ has length at least three. 
We will think of the $\Gamma_i$ as subsets of $C$.
We have:
\begin{lem} \label{lemcccy} Let $n=6$.
(i) If $\{i,j\} - \{k,m\} - \{r,s\}$ represent three consecutive vertices in a cycle $\Gamma_u$, then (up to a permutation) we have one of 

\noindent (a) $\{1,2\}-\{3,4\}-\{5,6\}$; 

\noindent (b)
$\{1,2\}-\{3,4\}-\{1,5\}$.

(ii) If $\{i,j\} - \{k,m\} - \{r,s\}-\{u,v\}$ represent four consecutive vertices in a cycle $\Gamma_u$, then (up to a permutation) we have one of

\noindent (a) $\{1,2\}-\{3,4\}-\{1,5\}-\{2,6\}$; 

\noindent (b)
$\{1,2\}-\{3,4\}-\{1,5\}-\{2,4\}$; 

\noindent(c)
$\{1,2\}-\{3,4\}-\{1,5\}-\{3,4\}$; 

\noindent (d)
$\{1,2\}-\{3,4\}-\{5,6\}-\{3,6\}$.

Further, in  a path   $\Gamma_u$  of length at least $4$ with  distinct  non-consecutive edges   $\{i,j\} - \{k,m\}$ and $\{i',j'\} - \{k',m'\}$,    the product 
$(i,j)(k,m)\cdot (i',j')(k',m')$ is not  in $\mathcal C_{(2^2,1^2)}$. 

(iii) If $\Gamma_i, \Gamma_j$ are distinct cycles, then no term of $\overline {\Gamma_i}\cdot \overline {\Gamma_j}$ is in $\mathcal C_{(2^2,1^2)}$.

(iv) If $ \Gamma_i$ has length greater than $4$, then the  number of elements of the form  $\alpha\beta, (\alpha,\beta \in \Gamma_i),$ that are in $\mathcal C_{(2^2,1^2)}$ is exactly  $|\Gamma_i|$, and none of these elements are in $C$. 
If $ \Gamma_i$ has length $4$, then the  number of elements of the form  $\alpha\beta, (\alpha,\beta \in \Gamma_i),$ that are in $\mathcal C_{(2^2,1^2)}$ is $2$, and none of these elements are in $C$.

%(v) For each $i \le r$, if $ \Gamma_i$ has length   $3$, then $\left (\overline {\Gamma_i}^2\right )_{\Gamma_i}=2\overline{\Gamma_i}.$
\end{lem}
\noindent {\it Proof}
(i) In such a path   $\{i,j\} - \{k,m\} - \{r,s\}$   we   have  $|\{r,s\} \cap \{i,j\}|=0,1$, which give the two cases.

(ii) The first part is easy to check using (i).
 If  all the vertices of   $\{i,j\} - \{k,m\},\{i',j'\} - \{k',m'\}$ are distinct, then one sees that $(i,j)(k,m)\cdot (i',j')(k',m')$ is not  in $\mathcal C_{(2^2,1^2)}$.

(iii) If $\Gamma_i, \Gamma_j$ are distinct cycles with $\alpha=\alpha_1\alpha_2\in \Gamma_1, \beta=\beta_1\beta_2 \in \Gamma_j, \alpha_k,\beta_k \in \mathcal C_2,k=1,2$, then 
$\{\alpha_1,\alpha_2\} \cap \{\beta_1,\beta_2\} = \emptyset$, since $\Gamma_i, \Gamma_j$ do not share a vertex.
Thus no element of $\overline{\Gamma_i}\cdot \overline{\Gamma_j}$ is in $\mathcal C_{(2^2,1^2)}$.

  (iv)  This now follows from (ii), which shows that only products of the form $(i,j)(r,s)\cdot (r,s)(u,v)$ will produce elements of $\mathcal C_{2,2}$. 
\qed\medskip

If we have a cycle $\Gamma_i$  of size greater than three, then 
 Lemma \ref {lemcccy}  (iv) shows that all cycles $\Gamma_j$ have size greater than $3$.
 The same result shows that no cycle has length four (since then we would have a principal set in $\mathcal C_{(2^2,1^2)}$ of size less than $15$.
  Thus $(\overline{C}^2)_{(2^2,1^2)}  =0\overline {C}+\overline {D},$
 where $|D|=15, C \cap D=\emptyset$. 
 We then similarly have  $(\overline{D}^2)_{(2^2,1^2)}  =0\overline {D}+\overline {C}.$ Thus $\mathcal C_{(2^2,1^2)}\setminus (C \cup D)$ has size $15$, and all of its cycles must have size $3$.
 Thus we may assume that  all cycles $\Gamma_i$ of $C$ have size $3$.
 
 In this case it is easy to see that the components $\Gamma_i$ of $C$ are as shown in Figure 3 (up to conjugacy):

   \begin{figure}[ht]
 \scalebox{.7}
{\includegraphics{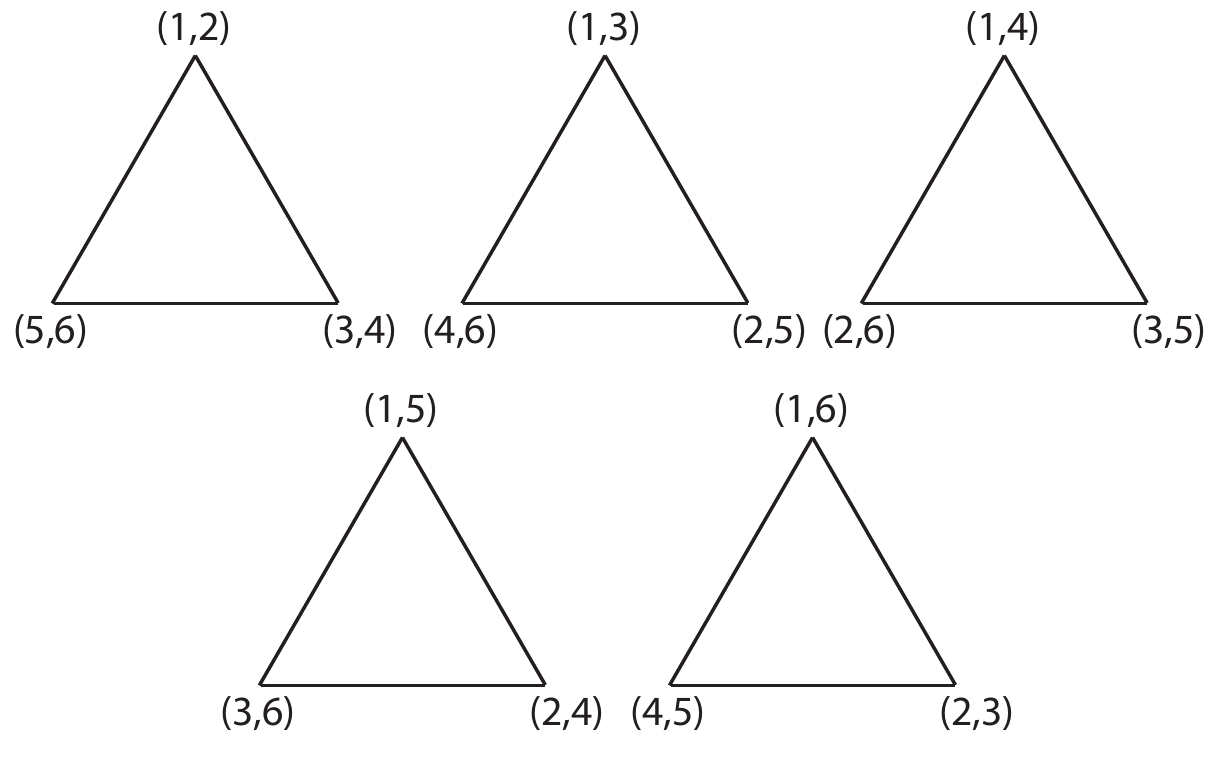}}
 \caption{  } 
 \end{figure}
 One now finds that the S-ring determined by $\overline {C}$ and $Z(\mathbb CS_6)$ is $\mathfrak S(S_6,H_{120})$, and that the complement of $C$ in $\mathcal C_{(2^2,1^2)}$ is also a principal set (of size $30$). Thus we cannot have any principal sets $C$ in $\mathfrak S$ with $|C|=15$ and with cycles $\Gamma_i$ of length greater than three.

This completes the proof of   Lemma \ref {lem333} and our discussion of the case $\lambda=2$.\qed\medskip

\medskip

For $n=6$ the element $\overline{\mathcal C_{(2^2,1^2)}}$ gives a $\lambda=6$ cover, thus if we have a cover $C$ with $\lambda=4$  then    
$\overline {\mathcal C_{(2^2,1^2)}}- \overline C$ has 
 $\lambda=2$,
and the S-ring generated by $Z(\mathbb CS_6)$ and $\overline C$ is the same as the S-ring generated by  $Z(\mathbb CS_6)$ and $\overline {\mathcal C_{(2^2,1^2)}}-\overline C$.
 Thus we have:

 \begin{prop} \label{prop2322} Any commutative S-ring $\mathfrak S$ over $S_6$ that contains $\mathcal C_2$ as a principal set, also contains $\mathcal C_3$ as a principal set, and either $\mathcal C_{2,2}$ 
  is a  principal set or $\mathfrak S$ contains a conjugate of $\mathfrak S(S_6,H_{120})$.
  
  Lastly, any principal set of $\mathfrak S$ that is properly contained in $\mathcal C_{2,2}$ is an orbit of one of its elements under the action of a conjugate of $H_{120}$. 
\qed \end{prop}

\noindent {\bf The $\mathcal C_{(2^3)}$ case}

\begin{lem} \label{lem222cyc}
The only non-empty and proper subsets $C$  of $\mathcal C_{(2^3)}\subset S_6$ that satisfy
$$\left ( \overline{C}\cdot \overline{\mathcal C_{2,2}}\right )_{(2,1^{4})}=\lambda  \overline{\mathcal C_2}
$$ for some $\lambda \ge 1$ are conjugates of
\begin{align*}
&(i) \,\,\, C_1=\{(1, 2)(3, 4)(5, 6), (1, 6)(2, 3)(4, 5),
(1, 5)(2, 4)(3, 6), (1, 4)(2, 6)(3, 5), \\& \qquad (1, 3)(2, 5)(4, 
    6)\}; or \\&
(ii) \,\,\, C_2=\mathcal C_{(2^3)} \setminus C_1.
\end{align*}
Further, any conjugate of   (i) or  (ii) generates (together with $Z(\mathbb CS_6)$) an  S-ring conjugate to  $\mathfrak S(S_6,H_{120})$.\end{lem}
\noindent {\it Proof} Since $|\mathcal C_{(2^3)}|=15$ we may assume that  $|C|\le 7$ and that $(1,2)(3,4)(5,6) \in C$.
For $\alpha \in \mathcal C_{(2^3)}$ one sees that $[\alpha \cdot \overline{\mathcal C_{(2^2,1^2)}}]_{(2,1^4)}$ is a sum of three elements. Thus
$3|C|=15\lambda$ i.e. $|C|=5\lambda$. Since $|C|\le 7$ we see that $\lambda=1$ and $|C|=5$. 
It is now easy to see  that up to a permutation we can assume that $(1,6)(2,3)(4,5) \in C$. Then there are two cases: (a) $(1,4)(2,6)(3,5) \in C$; and (b) $(1,4)(2,5)(3,6) \in C$. Either case 
quickly gives a unique $C$ of size $5$, these $C$s being conjugate.  One checks that $C_1$ and $C_2$ are principal sets of $\mathfrak S(S_6,H_{120})$.  This gives the result.
\qed\medskip 

This gives:

 \begin{prop} \label{prop2322} Let $\mathfrak S$ be a  commutative S-ring over $S_6$ that contains $\mathcal C_2$ as a principal set. Then either
  
  \noindent (i)  $\mathcal C_3, \mathcal C_{(2^3)}$ and $\mathcal C_{2,2}$ 
are   principal sets of $\mathfrak S$, or 
  
\noindent  (ii) $\mathfrak S$ contains $\mathfrak S(S_6,H_{120})$. 

 Lastly, any principal set of $\mathfrak S$ that is properly contained in $\mathcal C_{(2^3)}$ is an orbit of one of its elements under the action of a conjugate of $H_{120}$. 
\qed \end{prop}

\noindent {\bf The $\mathcal C_{(3,2,1)}$ case}

Let $\mathfrak S$ be a commutative S-ring over $S_n, n   \ge 5,$ that contains $\mathcal C_2, \mathcal C_3$ as  principal sets.
For  a principal set  $C \subseteq \mathcal C_{(3,2,1^{n-5})}\subset S_n$ we define sets $D_{i,j} \subset \mathcal C_3\subset S_n$ by 
$$\overline{C}= \sum_{1\le i< j \le n} (i,j)\overline {D_{i,j}}.
$$

Since $[(i,j)(r,s,t)\cdot \overline{\mathcal C_3}]_{(2,1^{n-2})}=(i,j)$,   and $\mathcal C_2$ is a principal set of $\mathfrak S$ we see that there is some $\lambda\in \mathbb N$ such that 
$$\left (\overline{C}\cdot \overline{\mathcal C_3}\right )_{(2,1^{n-2})}=\lambda \overline {\mathcal C_2}.$$ Further, from the definition of the $D_{i,j}$ we see that
$\lambda=|D_{i,j}|$ for all $1\le i<j\le n$. 

Since $[(i,j)(r,s,t)\cdot \overline {\mathcal C_2}]_{(3,1^{n-3})}=(r,s,t)$   and $\mathcal C_3$ is a principal set of $\mathfrak S$ we also have  $\mu \in \mathbb N$ such that 
$$\left (\overline{C}\cdot \overline{\mathcal C_2}\right )_{(3,1^{n-3})}=\mu \overline {\mathcal C_3}.$$
But we also have
$$\left (\overline{C}\cdot \overline{\mathcal C_2}\right )_{(3,1^{n-3})}= \sum_{1\le i<  j \le n} \overline {D_{i,j}}.
$$
Thus from $\mu \overline {\mathcal C_3}=\sum_{1\le i< j \le n} \overline {D_{i,j}}$ we see that for all $1\le i<j\le n$ we have 
$$\binom {n}{2} |D_{i,j}|=\binom {n}{2} \lambda=\mu |\mathcal C_3|=\mu\frac {n(n-1)(n-2)}{3}.$$
From this we obtain
\begin{lem}\label{lemsizs6} $\lambda=|D_{i,j}|$ and 
$$3\cdot |D_{i,j}|=2\mu(n-2).\qed 
$$ 
\end{lem}
Since any $(r,s,t) \in D_{i,j}$ must satisfy $i,j \notin \{r,s,t\}$ we see that
$$|D_{i,j}|\le \frac {(n-2)(n-3)(n-4)}{3}.$$
When $n=6$ this gives 
$|D_{i,j}|\le 8.$
But if $n=6$, then Lemma \ref {lemsizs6} shows that $3$ divides $\mu$ so that 
$$|D_{i,j}|=2\cdot \frac \mu 3 \cdot 4\ge 8.$$

It follows from these two equations that $|D_{i,j}|=8$ for all $1\le i<j\le n=6$, so that $D_{i,j}$ is maximal, i.e. $D_{i,j}=\{(r,s,t)^{\pm 1}:1\le r<s<t \le 6, i,j \notin \{r,s,t\}\}$,
and   $C = \mathcal C_{(3,2,1)}$.
%We note that $\mathcal C_{(3,2,1)}$ is a principal set of $\mathfrak S(S_6,H_{120})$.
Thus we now have shown 

 \begin{prop} \label{prop2322} Let $\mathfrak S$ be a  commutative S-ring over $S_6$ that contains $\mathcal C_2$ as a principal set. Then 
  $\mathcal C_{(3,2,1)}$ is a principal set of $\mathfrak S$, and 
 either
 
  \noindent (i)  $\mathcal C_3, \mathcal C_{(2^3)}, \mathcal C_{(3,2,1)}$ and $\mathcal C_{2,2}$ 
are   principal sets of $\mathfrak S$, or 
  
\noindent  (ii) $\mathfrak S$ contains a conjugate of $\mathfrak S(S_6,H_{120})$. 
\qed \end{prop}

    We note that in (ii) the only non-trivial conjugacy classes that do not split  are $\mathcal C_2,\mathcal C_3$ and $\mathcal C_{(3,2,1)}$.
    
\bigskip

\noindent {\bf The $\mathcal C_{(4,1^2)}$ case}

Let $C \subseteq \mathcal C_{(4,1^2)}$
   be a principal set.
Then  the fact that
$$\left ( (i,j,k,m)\cdot \overline {\mathcal C_3}\right )_{(2,1^4)}=(i,j)+(j,k)+(k,m)+(m,i),$$ shows that
$C$ determines a cover of $\lambda K_n$ by $4$-cycles. Thus, if $ [\overline C \cdot  \overline {\mathcal C_3}]_{(2,1^4)}=\lambda  \overline {\mathcal C_2}$, then
$$ 4\cdot |C|=15 \lambda.
$$

Thus  $15$ divides $|C|$.
Now  there are twelve  $\alpha \in \mathcal C_{(4,1^2)}$ such that $(1,2,3,4)\alpha$ is a $3$-cycle.
Thus considering  $\left ( \overline{C}\cdot \overline {\mathcal C_4}\right )_{(3,1^{3})}$   we see that (since $|\mathcal C_3|=40$) there is some $\mu \ge 1$ such that:
$$12\cdot |C|=40\mu,$$ (since $\mathcal C_3$ is a principal set) which shows that $10$ divides $|C|$. 

Thus $|C|$ is divisible by $30$. Note that $|\mathcal C_{(4,1^2)}|=90$.  If $|C|=60$, then we replace $C$ by $\mathcal C_{(4,1^2)} \setminus C$, so as to have $|C|=30$.

Now if $C^{-1} \ne C$, then $\mathcal C_{(4,1^2)} \setminus (C \cup C^{-1})$ is  a principal set (since it has size $30$) that is its own inverse. Thus we may assume that 
the principal set $C$ satisfies $|C|=30$ and  $C=C^{-1}$.

Write $\overline {C}=\alpha_1+\alpha_1^{-1}+\dots+\alpha_{15}+\alpha_{15}^{-1}$, so that 
 $(\overline{C}^2)_{(2^2,1^{2})}=2\alpha_1^2+2\alpha_2^{2}+\dots+2\alpha_{15}^2$, which shows that
  $\{\alpha_1^2,\alpha_2^{2},\dots,\alpha_{15}^2\} \subset \mathcal C_{(2^2,1^2)}$ is a principal set of $\mathfrak S$ of size $15$.
  Thus  Lemma  \ref {lem333} shows that such a set is conjugate to the set of $15$ elements given in that lemma. Call this set of $15$ elements $W_{15}$.

We note that the map 
$$\iota^2:\mathcal C_{(4,1^2)} \to \mathcal C_{(2^2,1^2)},\quad  \alpha \mapsto \alpha^2,$$ is a $2$ to $1$ surjection, where $(\iota^2)^{-1}((i,k)(j,m))=\{(i,j,k,m),(i,m,k,j)\}$, is an inverse pair.  
Since $C^{-1}=C$ we see that $C$ is completely determined by $\iota^2(C)$, which is conjugate to $W_{15}$. 
Thus $C$ must be a conjugate of
\begin{align} \tag {4.1} \{&
 (1, 2, 5, 4),(1, 4, 5, 2),
(1, 4, 3, 6),
    (2, 6, 5, 4),
        (1, 3, 4, 5),
    (1, 4, 6, 5),
    (2, 5, 3, 4),\\&
    (2, 3, 4, 6),
    (3, 5, 4, 6),\notag
    (1, 4, 2, 3),
    (1, 2, 4, 6),
    (1, 5, 4, 3),
    (3, 6, 4, 5),
    (2, 3, 6, 5),\\&
    (1, 6, 2, 5),
    (2, 4, 5, 6),
    (1, 3, 5, 6),\notag
    (2, 4, 3, 5),
    (1, 3, 6, 2),
    (1, 2, 3, 5),
    (1, 6, 3, 4),\\&
    (2, 5, 6, 3),
    (1, 5, 3, 2),
    (1, 5, 6, 4),\notag
    (1, 2, 6, 3),
    (1, 3, 2, 4),
    (2, 6, 4, 3),
    (1, 6, 4, 2),\\&
    (1, 5, 2, 6),
    (1, 6, 5, 3)\notag
\}.
\end{align}
Now one finds that the S-ring generated by $Z(\mathbb CS_6)$ and $\overline{C}$ is $\mathfrak S(S_6,H_{120})$. 
Using the fact that any conjugate of $\overline {W_{15}}$ that is not equal to $\overline {W_{15}}$ does not commute with $\overline {W_{15}}$, 
this now easily gives:
 \begin{prop} \label{prop23224} Let $\mathfrak S$ be a  commutative S-ring over $S_6$ that contains $\mathcal C_2$ as a principal set. Then either
  
  \noindent (i)  $\mathcal C_3, \mathcal C_{(2^3)}, \mathcal C_{(3,2,1)}, \mathcal C_{(4,1^2)}$ and $\mathcal C_{2,2}$ 
are   principal sets of $\mathfrak S$, or 
  
\noindent  (ii) $\mathfrak S$ contains a conjugate of $\mathfrak S(S_6,H_{120})$. 

Lastly, any principal set of $\mathfrak S$ that is properly contained in $\mathcal C_{(4^2,1^2)}$ is an orbit of one of its elements under the action of a conjugate of $H_{120}$, and is either a conjugate of (4.1) or of its complement in $\mathcal C_{(4,1^2)}$.
\qed \end{prop}

 \noindent {\bf The $\mathcal C_{(4,2)}$ case}
 
 There is a natural bijection $\mathcal C_{4,2} \to \mathcal  C_4, (i,j,k,m)(r,s) \mapsto (i,j,k,m)$. 
 So that the image of the  set of elements (4.1) of $\mathcal C_{(4,1^2)}$ is
 \begin{align} \tag {4.2} \{&
 (1, 2, 5, 4)(3,6),
 (1, 4, 5, 2)(3,6),
(1, 4, 3, 6)(2,5),
    (2, 6, 5, 4)(1,3),
        (1, 3, 4, 5)(2,6),\\&
    (1, 4, 6, 5)(2,3),
    (2, 5, 3, 4)(1,6),\notag
    (2, 3, 4, 6)(1,5),
    (3, 5, 4, 6)(1,2),
    (1, 4, 2, 3)(5,6),\\&
    (1, 2, 4, 6)(3,5),
    (1, 5, 4, 3)(2,6),\notag
    (3, 6, 4, 5)(1,2),
    (2, 3, 6, 5)(1,4),
    (1, 6, 2, 5)(3,4),\\&
    (2, 4, 5, 6)(1,3),
    (1, 3, 5, 6)(2,4),
    (2, 4, 3, 5)(1,6),\notag
    (1, 3, 6, 2)(4,5),
    (1, 2, 3, 5)(4,6),\\&
    (1, 6, 3, 4)(2,5),
    (2, 5, 6, 3)(1,4),
    (1, 5, 3, 2)(4,6),
    (1, 5, 6, 4)(2,3),\notag
    (1, 2, 6, 3)(4,5),\\&
    (1, 3, 2, 4)(5,6),
    (2, 6, 4, 3)(1,5),\notag
    (1, 6, 4, 2)(3,5),
    (1, 5, 2, 6)(3,4),
    (1, 6, 5, 3)(2,4)
\}.
\end{align}

 This case then follows by following the proof of the $\mathcal C_4$ case just completed.
This then gives:
 \begin{prop} \label{prop23225} Let $\mathfrak S$ be a  commutative S-ring over $S_6$ that contains $\mathcal C_2$ as a principal set. Then either
  
  \noindent (i)  $\mathcal C_3, \mathcal C_{(2^3)}, \mathcal C_{(3,2,1)}, \mathcal C_{(4,1^2)}, \mathcal C_{(4,2)}$ and $\mathcal C_{2,2}$ 
are   principal sets of $\mathfrak S$, or 
  
\noindent  (ii) $\mathfrak S$ contains $\mathfrak S(S_6,H_{120})$. 

Lastly, any principal set of $\mathfrak S$ that is properly contained in $\mathcal C_{(4^2,2)}$ is an orbit of one of its elements under the action of a conjugate of $H_{120}$, and is either a conjugate of (4.2) or of its complement in $\mathcal C_{(4,2)}$.
\qed \end{prop}

 \noindent {\bf The $\mathcal C_{(3^2)}$ case}

Consider a principal set $C \subseteq \mathcal C_{(3^2)}$. Note that if $\alpha=ab \in \mathcal C_{(3^2)}$ where $a,b \in \mathcal C_{3}$, and if $c \in \mathcal C_{3}$ with $\alpha \cdot c \in \mathcal C_{3}$, then $c \in \{a^{-1},b^{-1}\}$. 

Let $\overline C=\sum_{i=1}^r \alpha_i, \alpha_i=a_ib_i\in \mathcal C_{(3^2)}, a_i,b_i  \in \mathcal C_{3}$. Then
$$\left ( \overline{C}\cdot \overline {\mathcal C_{3}}\right)_{(3,1^{3})}=\sum_{i=1}^r (a_i+b_i).$$
Since $\mathcal C_{3}$ is a principal set we see that
$$\sum_{i=1}^r (a_i+b_i)=\lambda\overline{\mathcal C_3}.$$
Since $|\mathcal C_3|=40$ this gives 
$20\lambda=r.$ 
Now $r \le |\mathcal C_{(3^2)}|=40,$ so that we have $\lambda \le 2.$
If $\lambda=2$, then $r=40$ and $C = \mathcal C_{(3^2)}$. 

So assume that $\lambda=1, r=|C|=20$. We may also assume (by conjugating) that $(1,2,3)(4,5,6) \in C$.
We have  $[\overline{C}\cdot \overline {\mathcal C_3}]_{(3,1^3)}= 1\overline {\mathcal C_3}$, thus for each $\alpha \in\mathcal C_3$ there is a unique $\alpha^* \in
\mathcal C_3$ such that $\alpha \alpha^* \in C$. Since $(1,2,3)(4,5,6) \in C$ we  cannot have $(1,3,2)(4,5,6) \in C$, and so we must have $(1,3,2)(4,6,5) \in C$. Thus $C=C^{-1}$.

For $\alpha=ab \in \mathcal C_{(3^2)}, a,b \in \mathcal C_3,$ we let $S(\alpha)=\{ab,a^{-1}b,ab^{-1},a^{-1}b^{-1}\}$. 
Then the $S(\alpha), \alpha \in  \mathcal C_3,$ partition $\mathcal C_{(3^2)}$ into ten subsets $S_1,\dots,S_{10}$ of size $4$. 
Then any set $C \subset \mathcal C_{(3^2)}$ as above is a union of an inverse  pair $\{\alpha_i,\alpha_i^{-1}\} \subset S_i$
from each $S_i$. There are thus $2^{10}$ possible such $C$'s, but only $2^9$ if one assumes that (say) $(1,2,3)(4,5,6) \in C$

Of these $2^{9}$ choices for $C$ one finds that only twelve  satisfy
$\left ( \overline {C}^2\right )_{(3,3)}=\lambda \overline {C}+\mu \, \overline {\mathcal C_{(3^2)} \setminus C} $ for some $\lambda,\mu \ge 0$. 
These twelve  elements are in three conjugacy classes.
Of the representatives from these three conjugacy classes 
one checks that the subalgebra  of $\mathbb C S_6$ generated by such a 
$\overline C$ and $Z(\mathbb C S_6)$ is only an S-ring if (up to conjugacy) we have one of the following two cases.
% and one checks (\cite {ma}) that the subalgebra  of $\mathbb C S_6$ generated by such a 
%$\overline C$ and $Z(\mathbb C S_6)$ is only an S-ring if (up to conjugacy) we have
\begin{align}\tag {4.3}
\overline C_1&=  (1, 3, 5)(2, 4, 6) + (1, 5, 3)(2, 6, 4) + (1, 6, 4)(2, 5, 3) + (1, 4, 5)(2, 3, 6)\\& \notag + (1, 2, 6)(3, 5, 4) + (1, 5, 2)(3, 4, 6) + (1, 3, 4)(2, 5, 6) + (1, 6, 5)(2, 4, 3) \\& \notag
+ (1, 6, 3)(2, 4, 5) + (1, 3, 6)(2, 5, 4) + (1, 3, 2)(4, 6, 5) + 
    (1, 5, 6)(2, 3, 4)\\& \notag + (1, 4, 6)(2, 3, 5) + (1, 5, 4)(2, 6, 3) + (1, 2, 5)(3, 6, 4) + (1, 4, 2)(3, 6, 5) \\&+ (1, 2, 3)(4, 5, 6) + (1, 2, 4)(3, 5, 6) + (1, 4, 3)(2, 6, 5) + (1, 6, 2)(3, 4, 5); \notag \text { or }\\ \overline C_2&= \notag
     (1, 4, 5)(2, 3, 6) + (1, 2, 6)(3, 5, 4) +  \notag
    (1, 5, 2)(3, 4, 6) + (1, 6, 4)(2, 3, 5)\\& \notag + (1, 6, 5)(2, 3, 4) + (1, 6, 3)(2, 4, 5) + (1, 4, 3)(2, 5, 6) + (1, 3, 6)(2, 5, 4) \\& \notag+ (1, 5, 6)(2, 4, 3) + (1, 3, 2)(4, 6, 5) + (1, 5, 4)(2, 6, 3) + (1, 5, 3)(2, 4, 6)\\&  \notag + (1, 4, 6)(2, 5, 3) + 
(1, 2, 5)(3, 6, 4) + (1, 3, 4)(2, 6, 5) + (1, 4, 2)(3, 6, 5)\\&  \notag + (1, 2, 3)(4, 5, 6) + (1, 2, 4)(3, 5, 6) + (1, 3, 5)(2, 6, 4) + (1, 6, 2)(3, 4, 5).
\end{align}
However in each of these two cases the S-ring that one obtains is just a conjugate of $\mathfrak S(S_6,H_{120})$ again. 
\medskip

Thus we now  have:

 \begin{prop} \label{prop232251} Let $\mathfrak S$ be a  commutative S-ring over $S_6$ that contains $\mathcal C_2$ as a principal set. Then either
  
  \noindent (i)  $\mathcal C_3, \mathcal C_{(2^3)}, \mathcal C_{(3,2,1)}, \mathcal C_{(4,1^2)}, \mathcal C_{(4,2)}$, $\mathcal C_{2,2}$  and $\mathcal C_{(3^2)}$
are   principal sets of $\mathfrak S$, or 
  
\noindent  (ii) $\mathfrak S$ contains a conjugate of $\mathfrak S(S_6,H_{120})$. 

Lastly, any principal set of $\mathfrak S$ that is properly contained in $\mathcal C_{(3^2)}$ is an orbit of one of its elements under the action of a conjugate of $H_{120}$, and is either a conjugate of  one of the elements shown in (4.3)  or of its complement in $\mathcal C_{(3^2)}$.
\qed \end{prop}

   \medskip

 \noindent {\bf The $\mathcal C_{(5,1)}$ case}

Consider a principal set $C \subseteq \mathcal C_{5}=\mathcal C_{(5,1)}$. Now there are five elements  $\alpha \in \mathcal C_{2}$
such that $(1,2,3,4,5)\alpha \in \mathcal C_{(3,2,1)}$; since $ \mathcal C_{(3,2,1)}$ is a principal set we see that
$5\cdot |C|=\lambda \cdot | \mathcal C_{(3,2,1)}|=120\lambda$. This gives $|C|=24\lambda$.

Also, there are five   $\alpha \in \mathcal C_{3}$
such that $(1,2,3,4,5)\alpha \in \mathcal C_{(3,1^3)}$. 
Since $ \mathcal C_{(3,1^3)}$ is a principal set we see that
$5\cdot |C|=\mu \cdot | \mathcal C_{(3,1^3)}|=40\mu$.

\noindent {\bf Case 1:} $|C|=24$. Here $\lambda=1,\mu=3,$ and so we have (i) $[\overline {C}\cdot \overline{ \mathcal C_{2}}]_{(3,2,1)}=\overline{ \mathcal C_{(3,2,1)}}$ and (ii) 
$[\overline {C}\cdot \overline{ \mathcal C_{3}}]_{(3,1^3)}=3\overline{ \mathcal C_{(3,1^3)}}$.
 
 Let $\alpha_1,\dots,\alpha_{72} \in \mathcal C_{(5,1)}$ represent the inverse pairs in $\mathcal C_{(5,1)}$, so that $\mathcal C_{(5,1)}=\cup_{i=1}^{72} \{\alpha_i,\alpha_i^{-1}\}.$ 
 
 Assume first that $C^{-1}=C$. Then we can write $\overline {C}=\sum_{i=1}^{72} x_i(\alpha_i+\alpha_i^{-1}).$ Here $x_i =0,1$ satisfy
 
\noindent  (a) $x_i^2-x_i=0$ for $1\le i\le 72$; 

\noindent   (b) $\sum_{i=1}^{72} x_i=12$; 

\noindent  (c) $\left ( \sum_{i=1}^{72} x_i(\alpha_i+\alpha_i^{-1}) \cdot \overline {\mathcal C_2}\right )_{(3,2,1)}= \overline {\mathcal C_{(3,2,1)}}$;

\noindent  (d)  $\left ( \sum_{i=1}^{72} x_i(\alpha_i+\alpha_i^{-1}) \cdot \overline {\mathcal C_3}\right )_{(3,1^3)}=3 \overline {\mathcal C_{3}}$.
 
 If we now think of the $x_i$ as indeterminates in the polynomial ring $R=\mathbb Q[x_1,\dots,x_{72}]$, then each of (a)--(d) gives relations satisfied by the $x_i$. Let $I$ 
  be the ideal of $R$ generated by these relations. Finding a Gr\"obner basis for $I$ (\cite {ma}) we see that all but five of the $x_i$ are determined (in fact by a conjugacy we can assume that one of these five  is equal to zero). This leaves a small number of possibilities for $C$, and one checks that they are all equal to a conjugate of the $H_{120}$ orbit of some element of $\mathcal C_{(5,1)}$. Further each such element, together with the class sums, generates an S-ring that is just some conjugate of $\mathfrak S(S_6,H_{120})$. This does the case $|C|=24, C=C^{-1}$.\medskip
  
   So now assume   $|C|=24, C^{-1}\ne C$. Then we  write $\overline {C}=\sum_{i=1}^{72} x_{2i-1}\alpha_{i}+x_{2i}\alpha_{i}^{-1}.$ Here $x_i =0,1$ satisfy

\noindent  (a) $x_i^2-x_i=0$ for $1\le i\le 144$, and $x_{2i-1}x_{2i}=0$ for $1\le i\le 144$; 

\noindent   (b) $\sum_{i=1}^{144} x_i=24$; 

\noindent  (c) $\left ( \sum_{i=1}^{72} (x_{2i-1}\alpha_i+x_{2i}\alpha_i^{-1}) \cdot \overline {\mathcal C_2}\right )_{(3,2,1)}= \overline {\mathcal C_{(3,2,1)}}$;

\noindent  (d)  $\left (  \sum_{i=1}^{72} (x_{2i-1}\alpha_i+x_{2i}\alpha_i^{-1}) \cdot \overline {\mathcal C_3}\right )_{(3,1^3)}=3 \overline {\mathcal C_{3}}$.
 
Finding the ideal of $R=\mathbb Q[x_1,\dots,x_{144}]$ determined by these equations gives no solutions.
The above shows that if there is a $C \subset \mathcal C_{(5,1)}, |C|=24$, then the S-ring contains a conjugate of $\mathfrak S (S_6,H_{120})$, and a principal set of size $24$ in $\mathcal C_{(5,1)}$ is an orbit under the same conjugate of $H_{120}$.\medskip

\noindent {\bf Case 2:} $|C|=48$.  Here $\lambda=2,\mu=6,$ and  we aim to show that no such set exists. First note that if $C\ne C^{-1}$, then $D=\mathcal C_{(5,1)} \setminus (C \cup C^{-1})$ 
satisfies $D=D^{-1}, |D|=48,$ and $\overline {D} $ is in the S-ring. If $D$ is a union of disjoint sets of size $24$, then these sets are conjugates of
the sets of size $24$ determined in Case 1. One checks, however, that no two of these sets (there are six of them) commute, and so we see that $D$ must also be a principal set.
Thus we can assume that $C=C^{-1}$. 
Then we can write $\overline {C}=\sum_{i=1}^{72} x_i(\alpha_i+\alpha_i^{-1}).$ Here, as in Case 1, the   $x_i =0,1$ satisfy
 
\noindent  (a) $x_i^2-x_i=0$ for $1\le i\le 72$; 

\noindent   (b) $\sum_{i=1}^{72} x_i=24$; 

\noindent  (c) $\left (\sum_{i=1}^{72} x_i(\alpha_i+\alpha_i^{-1}) \cdot \overline {\mathcal C_2}\right )_{(3,2,1)}= 2\overline {\mathcal C_{(3,2,1)}}$;

\noindent  (d)  $\left (\sum_{i=1}^{72} x_i(\alpha_i+\alpha_i^{-1}) \cdot \overline {\mathcal C_3}\right )_{(3,1^3)}=6 \overline {\mathcal C_{3}}$.

Let $I$ 
  be the ideal of $R=\mathbb Q[x_1,\dots,x_{72}]$ generated by the polynomials determined by these relations. Finding a Gr\"obner basis for $I$ (\cite {ma}) we see that all but twelve of the $x_i$ are determined (again by a conjugacy we can assume that one of these twelve  is equal to zero).  Choosing values in $\{0,1\}$ for these eleven variables, one checks the $2^{11}$ cases, and one finds  that 
  there are only 
    $30$ (non-trivial) elements $C$ (any other case gives $I=R$),
and for each of these $30$ cases one finds  those such that $\left ( \overline{C}^2\right)_{(3,1^3)} $ has the form $m\overline {\mathcal C_3}, m\ge 0$.  There are $10$ of these, and one finds that they are all conjugate. Thus, if $C$ is one of these, then we check that the S-ring generated by $Z(\mathbb C S_6)$ and $\overline{C}$ is commutative and has dimension $34$, but does not contain $\mathcal C_2$ as a principal set. 
(Note: each of  these ten conjugate $C$s does give a commutative S-ring containing $\overline {\mathcal C_2}$, however $\mathcal C_2$ partitions as  $\{(1,2),(3,4),(5,6)\}$ and the complement; see Corollary  \ref {corthm1} (iii).) 
Thus there are no principal sets of size $48$ in $\mathcal C_{(5,1)}$ if $\mathcal C_2$ is a principal set.\medskip

  \noindent {\bf Case 3:} $|C|=72$.  Here $\lambda=3,\mu=9.$ 
  
  Assume first that $C=C^{-1}$. Then we can write $\overline {C}=\sum_{i=1}^{72} x_i(\alpha_i+\alpha_i^{-1}).$ Here, as in the above,  $x_i =0,1,$ satisfy
 
\noindent  (a) $x_i^2-x_i=0$ for $1\le i\le 72$; 

\noindent   (b) $\sum_{i=1}^{72} x_i=36$; 

\noindent  (c) $\left (\sum_{i=1}^{72} x_i(\alpha_i+\alpha_i^{-1}) \cdot \overline {\mathcal C_2}\right )_{(3,2,1)}= 3\overline {\mathcal C_{(3,2,1)}}$;

\noindent  (d)  $\left (\sum_{i=1}^{72} x_i(\alpha_i+\alpha_i^{-1}) \cdot \overline {\mathcal C_3}\right )_{(3,1^3)}=9 \overline {\mathcal C_{3}}$;

\noindent  (e)  $\left (\sum_{i=1}^{72} x_i(\alpha_i+\alpha_i^{-1}) \cdot \overline {\mathcal C_4}\right )_{(2,1^4)}=24 \overline {\mathcal C_{2}}$;

\noindent  (f) Let $\left (\overline{C}^2\right )_{(3,1^3)}= \sum_{\alpha \in \mathcal C_3} x_\alpha \alpha$, where $x_\alpha \in 
R=\mathbb Q[x_1,\dots,x_{72}]$. Then we have $x_\alpha=x_\beta$ for all $\alpha, \beta \in \mathcal C_3$, since $\mathcal C_3$ is a principal set. 

Here (e) is obtained in the same way as (c) and (d). Also, (f) is a consequence of the fact that $\mathcal C_3$ is a principal set of $\mathfrak S$. Constructing the ideal $I$ (including all $x_\alpha-x_\beta$ as in (f)) again one finds that there are $8$ variables that determine the rest ($7$ if one uses conjugacy to set one of them to be zero).  Choosing values in $\{0,1\}$ for these variables 
   gives $2^7$ cases to consider, of which all but eleven give $I=R$.  These eleven $C$s are in two conjugacy classes of sizes $1$ and $10$. The element $C_1$, representing  the single conjugacy class, is the orbit of some element of $\mathcal C_{(5,1)}$ under the action of $A_6$.
   Let $C_2$ be one of the ten conjugate elements. Then there is a subgroup $H_{36}$ of size $36$ such that the orbits of $H_{36}$ on $\mathcal C_{(5,1)}$ are $O_1,O_2,O_3,O_4$ where
   $|O_i|=36, i=1,2,3,4$. Further, these sets $O_i$ can be numbered so that
   $$C_1=O_1 \cup O_2,\quad C_2=O_1\cup O_3.$$
   A representative for $H_{36}$ up to conjugacy is
   $$H_{36}=\langle (1,6)(2,3,4,5),(1,2,5,4)(3,6)\rangle.$$
   We note that if one uses $O_1 \cup O_4$, then one obtains a commutative S-ring of dimension $26$ 
   that contains $\{(1, 3), (3, 5), (1,5), (2, 6), (4, 6), (2, 4)\}$ and its complement in $\mathcal C_2$, as  principal sets - see Corollary \ref {corthm1} (iii).

  One finds that each of the S-rings $\langle Z(\mathbb CS_6),\overline {C_1}\rangle, \langle Z(\mathbb CS_6),\overline {C_2}\rangle$ has dimension $12$ (recall that $Z(\mathbb CS_6)$ has dimension $11$) where $C_i, \mathcal C_{(5,1)} \setminus C_i,i=1,2, $ are principal sets of these S-rings. This concludes the case where $|C|=72, C=C^{-1}$.  
  We will let $C_{36}$ denote the element $C_2$ in what follows, and we will denote the S-ring $\langle Z(\mathbb CS_6),\overline {C_{36}}\rangle$ by $\mathfrak S_{36}$.
  
  \medskip
  
  Now if $C\ne C^{-1}$, then   we can write $\overline {C}=\sum_{i=1}^{72} x_{2i-1}\alpha_{i}+x_{2i}\alpha_{i}^{-1}.$ Here $x_i =0,1$ satisfy

  \noindent  (a) $x_i^2-x_i=0$ for $1\le i\le 144$ and $x_{2i-1}x_{2i}=0$ for $1\le i\le 72$; 

\noindent   (b) $\sum_{i=1}^{72} x_i=72$; 

\noindent  (c) $\left (\sum_{i=1}^{72} (x_{2i-1}\alpha_{i}+x_{2i}\alpha_{i}^{-1}) \cdot \overline {\mathcal C_2}\right )_{(3,2,1)}= 3\overline {\mathcal C_{(3,2,1)}}$;

\noindent  (d)  $\left (\sum_{i=1}^{72} (x_{2i-1}\alpha_{i}+x_{2i}\alpha_{i}^{-1}) \cdot \overline {\mathcal C_3}\right )_{(3,1^3)}=9 \overline {\mathcal C_{3}}$;

\noindent  (e)  $\left (\sum_{i=1}^{72} (x_{2i-1}\alpha_{i}+x_{2i}\alpha_{i}^{-1}) \cdot \overline {\mathcal C_4}\right )_{(2,1^4)}=24 \overline {\mathcal C_{2}}$;

\noindent  (f) Let $\left (\overline{C}^2\right )_{(3,1^3)}= \sum_{\alpha \in \mathcal C_3} x_\alpha \alpha$, where $x_\alpha \in 
R=\mathbb Q[x_1,\dots,x_{72}]$. Then we have $x_\alpha=x_\beta$ for all $\alpha, \beta \in \mathcal C_3$, since $\mathcal C_3$ is a principal set. 

One then finds that the ideal determined by (a)--(f) is the whole ring, so that there are no solutions in this situation.

Now if $C\subsetneq \mathcal C_{(5,1)}$ is a principal set of $\mathfrak S$ with $|C| > 72$, then there are certainly principal sets $D$ of $\mathfrak S$
with $|D|<72$. However in the above considerations of Cases 1,2 we have shown that the S-ring generated by $D$ and $Z(\mathbb CS_6)$ does not have a principal set of size greater than $72$. Thus the situation $|C|>72$  does not occur, and we have now proved:

\begin{prop} \label{propC5} Let $\mathfrak S$ be a commutative S-ring over $S_6$ containing $\overline {\mathcal C_2}$. Then any principal set $C \subset \mathcal C_{(5,1)}$ of $\mathfrak S$  is either a conjugate of $C_{36}$, or is the orbit of one of its elements under the action of $A_6$ or of a conjugate of  $H_{120}$.\qed 
\end{prop}

  \medskip 
   \noindent {\bf The $\mathcal C_{(6)}$ case}

Consider a principal set $C \subsetneq \mathcal C_{(6)}$. Now there are six elements  $\alpha \in \mathcal C_{(4,1^2)}$
such that $(1,2,3,4,5,6)\alpha \in \mathcal C_{(3,1^3)}$; since $ \mathcal C_{(3,1^3)}$ is a principal set we see that
$6\cdot |C|=\lambda_1 \cdot | \mathcal C_{(3,1^3)}|=40\lambda_1$. This gives (i) $3|C|=20\lambda_1$.
  
  Similarly, by considering 
  
  \noindent (ii) $\left (\overline{C} \cdot \overline {\mathcal C_{(5,1)}}\right )_{(2,1^4)}=\lambda_2\overline {\mathcal C_2}$, we see that $2|C|=5\lambda_2$; 
  
   \noindent (iii) $\left (\overline{C} \cdot \overline {\mathcal C_{(3,2,1)}}\right )_{(3,1^3)}=\lambda_3\overline {\mathcal C_3}$, we see that $3|C|=10\lambda_3$; 
   
    \noindent (iv) $\left (\overline{C} \cdot \overline {\mathcal C_{2,2}}\right )_{(3,2,1)}=\lambda_4\overline {\mathcal C_{(3,2,1)}}$, we see that $3|C|=20\lambda_4$. 
   
  One sees that $|C|$ is divisible by $20$.  Also $|\mathcal C_{(6)}|=120$. Let $\alpha_i, i=1,\dots,60$ be representatives for the inverse pair sets in $\mathcal C_{(6)}$.

    \noindent {\bf Case 1:} $|C|=20$.  Here $\lambda_1=3,\lambda_2=8,\lambda_3=6,\lambda_4=3$. 
  
  First assume that $C=C^{-1}$. Then we can write $\overline {C}=\sum_{i=1}^{60} x_i(\alpha_i+\alpha_i^{-1}).$ Here, as in the above,  $x_i =0,1$ satisfy

\noindent  (a) $x_i^2-x_i=0$ for $1\le i\le 60$; 

\noindent   (b) $\sum_{i=1}^{60} x_i=10$; 
  
  \noindent   (c) relations for each of (i)--(iv) above;
  
  \noindent   (d) let $\left (\overline{C}^2\right )_{(3,1^3)}= \sum_{\alpha \in \mathcal C_3} x_\alpha \alpha$, where $x_\alpha \in 
R=\mathbb Q[x_1,\dots,x_{10}]$. Then we have $x_\alpha=x_\beta$ for all $\alpha, \beta \in \mathcal C_3$, since $\mathcal C_3$ is a principal set.
  
    \noindent   (e) Let $\left (\overline{C}^2\cdot \overline {\mathcal C_{2,2}}\right )_{(3,1^3)}= \sum_{\alpha \in \mathcal C_3} x_\alpha \alpha$, where $x_\alpha \in 
R=\mathbb Q[x_1,\dots,x_{10}]$. Then we have $x_\alpha=x_\beta$ for all $\alpha, \beta \in \mathcal C_3$, since $\mathcal C_3$ is a principal set.
  
Constructing the ideal $I$ (including all $x_\alpha-x_\beta$ as in (d),(e)) again one finds that there are $5$ variables that determine the rest ($4$
 if one uses conjugacy to set one of them to be zero). Looking at the $2^4$ cases one finds that there are $5$ possibilities for $C$. 
 These are all in a single conjugacy class, so we need only consider one of them, $C$ say. 
 One finds that the S-ring generated by $Z(\mathbb CS_6)$ and $\overline {C}$ is a conjugate of $\mathfrak S(S_6,H_{120})$. This does the case where $C=C^{-1}$.\medskip
 
 Now if $C\ne C^{-1}$, then   we can write $\overline {C}=\sum_{i=1}^{60} x_{2i-1}\alpha_{i}+x_{2i}\alpha_{i}^{-1}.$ 
Using the same tests as in the $C=C^{-1}$ case above one creates an ideal that has nine  variables that determine the rest. Enumerating the various possibilities for $C$ gives ten 
non-conjugate elements. One checks that each of these generates $\mathbb CS_6$ as an S-ring. This shows that this case cannot occur.
We further note
\begin{lem}\label{lemC620} If $C_1, C_2 \subseteq \mathcal C_6$ are distinct principal sets of a commutative S-ring over 
$S_6$ that contains $\mathcal C_2$ as a principal set, then at most one of $C_1, C_2$  can have size $20$.
\end{lem}
\noindent {\it Proof} If $|C_1|=20$, then the above shows that $C_2$ is a conjugate of $C_1$; but there are only six such conjugates and it is easy to check that no two such distinct commute.\qed

\medskip

 \noindent {\bf Case 2:} $|C|=40$.  Here $\lambda_1=6,\lambda_2=16,\lambda_3=12,\lambda_4=6$. 

  First assume that $C=C^{-1}$. Then we can write $\overline {C}=\sum_{i=1}^{60} x_i(\alpha_i+\alpha_i^{-1}).$ Applying the same ideal calculation as in the $|C|=20$ case 
  one obtains an ideal  having $4$ variables that determine the rest. 
Enumerating the various possibilities for $C$ gives four elements, all of them conjugate to each other.
Considering one of them, $C$ say, one  finds that the  S-ring generated by $Z(\mathbb CS_6)$ and $\overline {C}$ is a conjugate of $\mathfrak S(S_6,H_{120})$. This does the case where $C=C^{-1}$.\medskip

Repeating the above for the situation where $|C|=40, C \ne C^{-1},$ one finds that there are no solutions.\medskip 
  
   \noindent {\bf Case 3:} $|C|=60$.  Here $\lambda_1=9,\lambda_2=24,\lambda_3=18,\lambda_4=9$.
   
  First assume that $C=C^{-1}$. Then we can write $\overline {C}=\sum_{i=1}^{60} x_i(\alpha_i+\alpha_i^{-1}).$ Applying the same ideal calculation as in the $|C|=20$ case 
  one obtains an ideal  having $5$ variables that determine the rest. 
Enumerating the various possibilities for $C$ gives six possibilities for $C$, and there are two  conjugacy classes of such elements.
If  $C$ say, represents either of these classes, one  finds that the  S-ring generated by $Z(\mathbb CS_6)$ and $\overline {C}$ is a conjugate of $\mathfrak S(S_6,H_{120})$. This does the case where $C=C^{-1}$.  \medskip
   
Repeating the above for the situation where $|C|=60, C \ne C^{-1},$ one finds that there are no solutions.

The cases $|C|>60$ are dealt with as in the situation $|C|>72$ of Case 2.

\medskip 

This concludes consideration of all conjugacy classes of $S_6$ where $\mathfrak S$ has $\mathcal C_2$ as a principal set.

The following is mostly a summary of what we have done:

%%  Let $D_1,\dots,D_{19}$ be the classes of $\mathfrak S(S_6,H_{120})$.
 % Let $E_1,\dots,E_{11}$ be the classes of $\mathfrak S(S_6,A_6)$, where  $E_{10}, E_{11} \subset \mathcal C_{(5,1)}, |E_{10}|=|E_{11}|=72$.
   % Let $F_1,\dots,F_{11}$ be the classes of the other S-ring, where  $F_{10}, F_{11} \subset \mathcal C_{(5,1)}, |F{10}|=|F_{11}|=72$.

\begin{prop} \label{prop22} 
  Let $\mathfrak S$ be a commutative S-ring over $S_6$ containing $\mathcal C_2$ as a principal set. Then
  
  \noindent (i)  $\mathcal C_{(3,1^3)}  $ and  $\mathcal C_{(3,2,1)}$ are  principal sets of $\mathfrak S$.
 \smallskip

   \noindent (ii) If $C \subsetneq \mathcal C_{2,2}$ is a principal set of $\mathfrak S$, then $|C|=15$ or $|C|=30$, and $\overline {C}$ is either a conjugate of the element shown in 
   Lemma \ref{lem333}, or is   $\overline {\mathcal C_{2,2}}-\overline{C}$ for such a set of size $15$. There are six conjugates of each such set $C$.  No two distinct such conjugates commute.
 The S-ring generated by $Z(\mathbb CS_6)$ and $\overline {C}$ is $\mathfrak S(S_6,H_{120})$. 
   
  \smallskip

 \noindent (iii) If $C \subsetneq \mathcal C_{(2^3)}$ is a principal set of $\mathfrak S$, then $|C|=5$ or $|C|=10$. If $|C|=5$, then $C$ is conjugate to the element shown in Lemma \ref {lem222cyc} (i).  There are six conjugates of each such set $C$.  No two distinct such conjugates commute.
 The S-ring generated by $Z(\mathbb CS_6)$ and $\overline {C}$ is  a conjugate of   $\mathfrak S(S_6,H_{120})$. 
   \smallskip \smallskip

 \noindent (iv) If $C \subsetneq \mathcal C_{(4,1^2)}$ is a principal set of $\mathfrak S$, then $|C|=30$ or $|C|=60$. If $|C|=30$, then $C$ is a conjugate of the element shown in (4.1), otherwise it is the complement in $\mathcal C_{(4,1^2)}$ of such an element.  There are six conjugates of each such set $C$. No two distinct such conjugates commute.
 The S-ring generated by $Z(\mathbb CS_6)$ and $\overline {C}$ is   a conjugate of  $\mathfrak S(S_6,H_{120})$. 
 \smallskip \smallskip

   \noindent (v) If $C \subsetneq \mathcal C_{(4,2)}$ is a principal set of $\mathfrak S$, then $|C|=30$ or $|C|=60$. 
    If $|C|=30$, then $C$ is a conjugate of the element shown in (4.2), otherwise it is the complement in $\mathcal C_{(4,2)}$ of such an element. There are six conjugates of each such set $C$. No two distinct such conjugates commute.
 The S-ring generated by $Z(\mathbb CS_6)$ and $\overline {C}$ is   a conjugate of  $\mathfrak S(S_6,H_{120})$. 
  \smallskip \smallskip

    \noindent (vi) If $C \subsetneq \mathcal C_{(3^2)}$ is a principal set of $\mathfrak S$, then $|C|=20$ and $\overline {C}$ is one of the two elements  $\overline {C_1}$, $\overline {C_2}$ shown in (4.3).  Now $\overline {C_1}$ and  $\overline {C_2}$ commute. Each of $C_1, C_2$ has six conjugates; no two distinct conjugates of each $\overline {C_i},i=1,2,$ commute, and each $\overline {C_1}^g$ only commutes with a single conjugate of  $\overline {C_2}$. 
  For $i=1,2$, the S-ring generated by $Z(\mathbb CS_6)$ and $\overline {C_i}$ is   a conjugate of  $\mathfrak S(S_6,H_{120})$.
   \smallskip \smallskip

   \noindent (vii) If $C \subsetneq \mathcal C_{(5,1)}$ is a principal set of $\mathfrak S$, then either 
   
   (viia) $|C|=24$ and $C$ is an orbit of an element of $\mathcal C_{(5,1)}$ under the action of a conjugate of $H_{120}$; or
   
   (viib)  $|C|=120$ and $C$ is an orbit of an element of $\mathcal C_{(5,1)}$ under the action of a conjugate of $H_{120}$; or 
   
(viic)   $|C|=72$ and $C$ is an $A_6$-orbit of an element  of $\mathcal C_{(5,1)}$.
% under the action of $A_6$; or 

(viid)   $|C|=72$ and $C$ is a conjugate of $C_{36}$.

Further, no two distinct conjugates of $C$, for $C$ of  type (viia) commute; no two distinct conjugates of  type (viib) commute;
no  conjugate of  type (viia) commutes with a conjugate of type (viib) unless their sum is $\overline {\mathcal C_{(5,1)}}$. Any conjugate of type (viic) commutes with any conjugate of type (viia) or (viib). 

Any element of type (viia) or (viib) generates (with $Z(\mathbb CS_6)$) an S-ring which is a conjugate of  $\mathfrak S(S_6,H_{120})$.
Any element of type  (viic) generates (with $Z(\mathbb CS_6)$) an S-ring which is a conjugate of  $\mathfrak S(S_6,A_6)$.
 \smallskip \smallskip

  \noindent (viii) If $C \subsetneq \mathcal C_{(6)}$ is a principal set of $\mathfrak S$, then  $|C|\in \{20,40,60\}$, and in each case  the S-ring generated by $Z(\mathbb CS_6)$ and $\overline {C}$ is   a conjugate of   $\mathfrak S(S_6,H_{120})$. Let $O_{20},O_{40},O_{60}$ denote the sets of conjugates (of the sums $\overline{C}$) for  each case, so that each $O_i$ has size six. Then for $i \in \{20,40,60\}$ no two distinct elements of $O_i$ commute, and for
  distinct $i,j \in \{20,40,60\}$ elements of $O_i, O_j$ commute if and only if the corresponding sets are disjoint.

   Further, the only principal elements of $\mathfrak S(S_6,H_{120})$ that $\overline {C_{36}}$ commutes with are $\overline {\mathcal C_2}, \overline {\mathcal C_3}$ and $\overline {\mathcal C_{(3,2,1)}}$. 
   
   %xxx
   
  \end{prop}
  \noindent {\it Proof} Here (i) follows from Proposition \ref{prop23} and Proposition \ref {prop2322}. Part (ii) follows from Lemma \ref {lem333},  the fact that $[S_6:H_{120}]=6$, and a calculation to show that distinct conjugates don't commute.
  (iii) follows from Lemma \ref {lem222cyc}, the fact that $[S_6:H_{120}]=6$, and a calculation to show that distinct conjugates don't commute.
  
  Part (iv) follows from the proof of Proposition  \ref {prop23224}, and a calculation to show that distinct conjugates don't commute.
  Part (v) follows from the discussion of the $\mathcal C_{(4,2)}$ case,   and a calculation to show that distinct conjugates don't commute.
   Part (vi) follows from the discussion of the $\mathcal C_{(3^2)}$ class that resulted in the two cases shown in (4.3) (up to conjugacy), and a calculation to show that distinct conjugates don't commute.
  Parts (vii) and (viii)  follow  from the discussion of the $\mathcal C_{(5,1)}$ and $\mathcal C_{(6)}$ classes, and a calculation to show that certain conjugates of these elements don't commute.\qed
  
  We use this to prove:

\begin{thm}\label{thmc2pset} The only commutative Schur rings over $S_6$ containing $ {\mathcal C_2}$ as a principal set are:

(1) $\mathfrak S(S_6,S_6)$;

(2) $\mathfrak S(S_6,H_{120})$;

(3) $\mathfrak S(S_6,A_6)$;

 (4)  $\mathfrak S_{36}$.
 
 \end{thm}
\noindent {\it Proof} The last statement of Proposition \ref {prop22} shows that such an S-ring cannot have split $S_6$ classes that are orbits under a conjugate of $H_{120}$, and also have $C_{36}$ (or its complement) in it. One checks that the S-ring generated by any conjugate of $\mathfrak S_{36}$ and      $\mathfrak S(S_6,A_{6})$ is not commutative. Given the fact that, according to Proposition \ref {prop22}, not very many 
  principal elements of conjugates of $\mathfrak S(S_6,H_{120})$
commute, it is now easy to prove the theorem.\qed

  \section {When $\mathcal C_2$ is not a principal set}

  We consider each case as enumerated in Corollary \ref {corthm1}.

  \noindent {\bf Case (ii):}  $\mathcal C_2=C_1 \cup C_2,$ where 
 $C_1=\{(1,2)\} \cup \{(i,j):3\le i<j\le 6\}$ and $
C_2=\mathcal C_2 \setminus C_1.$
  
  One finds that the S-ring generated by $Z(\mathbb CS_6)$, $\overline {C_1}$ and  $\overline {C_2}$ is $\mathfrak S(S_6,S_2 \times S_4)$, where 
  $S_2 \times S_4$ is the subgroup $\langle (1,2),(3,4),(3,4,5,6)\rangle\le S_6$. There are $34$ principal sets $O_1,\dots,O_{34}$ of sizes $1,3,6,8,12,16,24,48$. 
  The idea is to show that no proper, non-empty subset of each $O_i$ can be a principal set of a commutative S-ring containing $\mathfrak S(S_6,S_2 \times S_4)$.
  
  One can check this directly (using \cite {ma}) for each $O_i$ with $|O_i|\le 8$. For the rest let $O_i=\{o_1,\dots,o_m\}$. Let $R=\mathbb Q[x_1,\dots,x_m]$ be a polynomial ring, and let
  $e=\sum_{j=1}^m x_jo_j.$ One considers the ideal $I$ of $R$ generated by all the coefficients of the elements $e\overline{O_k}-\overline{O_k}e$, these being linear polynomials in the
  $x_j$. One finds that a Gr\"obner basis for $I$ has  only elements of the form $P_h-P_{u_k}$, where the $P_{u_k}$ are free variables, and there are at most six of the $P_{u_k}$. 
  One chooses a subset of the $P_{u_k}$, and puts these equal to $1$, while one puts the rest equal to $0$. This produces an ideal, that determines an element in $\mathbb CS_6$. One then sees whether this element (together with $Z(\mathbb CS_6))$ generates a commutative S-ring. One finds that  $\mathfrak S(S_6,S_2 \times S_4)$ is the only possible case.\medskip

  \noindent {\bf Case (iii):}  $\mathcal C_2=C_1 \cup C_2$ where 
 $C_1=\{(1,2),(3,4),(5,6)\} $ and
$C_2=\mathcal C_2 \setminus C_1.$
  
   One finds that the S-ring generated by $Z(\mathbb CS_6)$, $\overline {C_1}$ and  $\overline {C_2}$ is $\mathfrak S(S_6,S_2\wr S_3)$, where 
  $S_2 \wr S_3$ is the (wreath product) subgroup $\langle (1,2),(3,4),(5,6),(1,3)(2,4),(3,5)(4,6)\rangle\le S_6$. There are $34$ principal sets $O_1,\dots,O_{34}$ of sizes $1,3,6,8,12,16,24,48$. 
  
    We note that $S_6$ has an outer automorphism $\alpha$ such that  $\alpha(\mathfrak S(S_6,S_2 \times S_4))=\mathfrak S(S_6,S_2\wr S_3)$. Thus this case  follows from Case (ii).

  %The idea is to show that no proper, non-empty subset of each $O_i$ can be a principal set of a commutative S-ring containing $\mathfrak S(S_6,S_2 \wr S_3)$.
  
  %One can checks this directly for each $O_i$ with $|O_i|\le 8$, and the rest of the cases are dealt with as in Case (ii) above. One finds that  $\mathfrak S(S_6,S_2 \wr S_3)$ is the only possible case.
  
  \medskip
  
   \noindent {\bf Case (iv):} $\mathcal C_2=C_1 \cup C_2$ where 
 $C_1=\{(1,2),(1,3),(2,3),(4,5),(4,6),(5,6)\}$ and $
C_2=\mathcal C_2 \setminus C_1.$  

   One finds that the S-ring generated by $Z(\mathbb CS_6)$, $\overline {C_1}$ and  $\overline {C_2}$ is $\mathfrak S(S_6,S_3\wr S_2)$, where 
  $S_3 \wr S_2$ is the (wreath product) subgroup $\langle (1,2),(2,3),(5,6),(1,4)(2,5)(3,6)\rangle\le S_6$. There are $26$ principal sets $O_1,\dots,O_{26}$ of sizes $1,4,6,9,12,18,36,72$. 
  The idea is to again show that no proper, non-empty subset of each $O_i$ can be a principal set of a commutative S-ring containing $\mathfrak S(S_6,S_3 \wr S_2)$.\medskip

     \noindent {\bf Case (v):} $\mathcal C_2=C_1 \cup C_2$ where 
 $C_1=\{(1,6),(2,6),(3,6),(4,6),(5,6)\}$ and $
C_2=\mathcal C_2 \setminus C_1.$  

   One finds that the S-ring generated by $Z(\mathbb CS_6)$, $\overline {C_1}$ and  $\overline {C_2}$ is $\mathfrak S(S_6,S_5)$. There are $19$ principal sets $O_1,\dots,O_{19}$ of sizes $1,5,10,15,20,24,30,40,60,120$. 
   We note that $S_6$ has an outer automorphism, and that  $\mathfrak S(S_6,H_{120})$ and  $\mathfrak S(S_6,S_5)$ are related by this automorphism. Thus this case follows from (iv). 
   
   This concludes the proof of Theorem \ref {lem2222}. 
  \qed \medskip 

  \noindent{\it Proof of Corollary  \ref {corS66}}  In \cite {hum2} it is shown that any commutative S-ring $\mathfrak S$, of maximal dimension $s_G$, over a group $G$  contains $Z(\mathbb CG)$. Thus the principal sets of $\mathfrak S$ give a partition of $G$ that is a refinement of the partition of $G$ by conjugacy classes. In particular, for $G=S_6$, we must have $\overline {\mathcal C_2} \in \mathfrak S$. Thus any such S-ring must be in the list given in Theorem \ref {lem2222}. However none of the S-rings listed in Theorem \ref {lem2222} has dimension $76=s_{S_6}$. This proves Corollary  \ref {corS66}.
  \qed

\end{document}